\documentclass[12pt,a4paper]{article}
\usepackage{amssymb}
\usepackage{amsfonts}

\usepackage{amsmath}
\usepackage{amsthm}
\usepackage{enumerate}

\RequirePackage[OT1]{fontenc}
\RequirePackage{amsthm,amsmath, xcolor}
\RequirePackage[numbers]{natbib}
\RequirePackage[colorlinks,citecolor=blue,urlcolor=blue]{hyperref}
\usepackage{amscd,amsfonts,amssymb,latexsym,array,hhline,graphicx}
\usepackage{float}

\newtheorem{theorem}{Theorem}[section]
\newtheorem{proposition}[theorem]{Proposition}

\newtheorem{lemma}[theorem]{Lemma}

\newtheorem{remark}{Remark}[section]

\newcommand\fdem{$\Box$}
\newcommand\cA{{\cal A}}

\newcommand\cP{{\cal P}}
\newcommand\cG{{\cal G}}
\newcommand\cF{{\cal F}}

\newcommand\cB{{\cal B}}

\newcommand\cR{{\cal R}}

\newcommand\e{\epsilon}
\newcommand\ve{\varepsilon}

\newcommand\Er{\mbox{Err}}

\def\bbr{{\mathbb R}}

\def\text#1{\hbox{#1}}
\def\proof{{\noindent \bf Proof. }}
\def\endproof{\mbox{\ $\qed$}}

\def\E{{\bf E}}

\def\P{{\bf P}}
\def\B{{\bf B}}

\def\l{{\bf l}}
\def\C{{\bf C}}
\def\D{{\bf D}}

\def\L{{\bf L}}

\def\m{{\bf m}}

\def\q{{\bf q}}
\def\u{{\bf u}}
\def\l{{\bf l}}
\def\b{{\bf b}}

\def\s{{\bf s}}

\def\Chi{{\bf 1}}

\def\d{\mathrm{d}}
\def\build #1_#2{\mathrel{\mathop{\kern 0pt #1}\limits_{#2}}}

\newcommand\Trg{\mbox{Tr}}
\newcommand{\wh}{\widehat}
\newcommand{\wt}{\widetilde}
\newcommand{\zs}[1]{{\mathchoice{#1}{#1}{\lower.25ex\hbox{$\scriptstyle#1$}}
{\lower0.25ex\hbox{$\scriptscriptstyle#1$}}}}

\numberwithin{equation}{section}

\begin{document}
\title{
Sequential Model Selection Method for Nonparametric Autoregression
\thanks{
This work was supported by
 the RSF grant 17-11-01049 (National Research Tomsk State University).
}
}

\author{Ouerdia Arkoun
\thanks{
Sup'Biotech, Laboratoire BIRL, 66 Rue Guy Moquet, 94800 Villejuif, France  and
Laboratoire de Math\'ematiques Raphael Salem, Normandie Universit\'e,
 UMR 6085 CNRS- Universit\'e de Rouen,  France,
  e-mail: ouerdia.arkoun@gmail.com
}
,
Jean-Yves Brua
\thanks{
Laboratoire de Math\'ematiques Raphael Salem, Normandie Universit\'e,
 UMR 6085 CNRS- Universit\'e de Rouen,  France,
  e-mail :  jean-yves.brua@univ-rouen.fr
}
 and\\
 Serguei Pergamenshchikov\thanks{
 Laboratoire de Math\'ematiques Raphael Salem,
 UMR 6085 CNRS- Universit\'e de Rouen Normandie,  France
and
International Laboratory of Statistics of Stochastic Processes and
Quantitative Finance, National Research Tomsk State University,
 e-mail:
Serge.Pergamenshchikov@univ-rouen.fr } 
}

\date{}

\maketitle

\begin{abstract}
In this paper for the first time the nonparametric autoregression estimation problem for the quadratic risks is considered.
To this end
we develop
 a new adaptive sequential model selection method
based 
on
 the efficient
sequential  kernel estimators proposed
 by Arkoun and Pergamenshchikov (2016).
 Moreover,
 we develop a new  analytical tool for general regression models
  to obtain
 the non asymptotic sharp oracle
 inequalities for both usual quadratic and robust 
 quadratic risks.  Then, we show that the constructed sequential model selection procedure is optimal in the sense of oracle
 inequalities.
 \end{abstract}

\vspace*{5mm}

\vspace*{5mm}

\newpage

\section{Introduction}

One of the standard linear models in general theory of time series is the
autoregressive model (see, for example, \cite{Anderson1970} and the references therein).
Natural extensions for such models are nonparametric autoregressive models which are defined
by
\begin{equation}\label{sec:In.1}
y_\zs{k} = S(x_\zs{k}) y_\zs{k-1}+ \xi_\zs{k}
\quad\mbox{and}\quad
x_\zs{k}=a+\frac{k(b-a)}{n}\,,
\end{equation}
where $S(\cdot)\in\L_\zs{2}[a,b]$ is unknown function,
$a<b$ are fixed known constants,
$1\le k\le n$,
 the initial value
  $y_\zs{0}$ is a constant and
 the noise $(\xi_\zs{k})_\zs{k\ge 1}$ is i.i.d. sequence of 
unobservable random variables with $\E\xi_\zs{1}=0$ and 
$\E\xi^{2}_\zs{1}=1$. 

The problem is to estimate the function $S(\cdot)$
on the basis of the observations $(y_\zs{k})_\zs{1\le k\le n}$ under the 
condition that the noise  distribution is unknown.

It should be noted that the varying coefficient
principle is well known in the regression analysis.
It permits the use of more complex forms for  regression coefficients
and, therefore, the models constructed via this method
are more adequate for applications (see, for example,
\cite{FanZhang2008}, \cite{LuYangZhou2009}). In this paper we consider
the varying coefficient autoregressive models \eqref{sec:In.1}.
There is a number of papers which consider these  models
such as  \cite{Dalhaus1996a}, \cite{Dalhaus1996b} and \cite{Belitser2000}.
In all these papers, the authors propose some asymptotic (as $n\to\infty$) methods  for different
identification studies without considering optimal estimation issues.
To our knowledge, for the first time, the minimax estimation
problem
for the model \eqref{sec:In.1}
has been treated in 
\cite{ArkounPergamenshchikov2008} and \cite{MoulinesPriouretRoueeff2005} in the nonadaptive case,
i.e. for the known regularity of the function
$S$.
Then, in  \cite{Arkoun2011} it is proposed to use the
sequential analysis method
for the adaptive  pointwise estimation  problem
in the case where the unknown H\"older regularity is less than one, i.e when the function $S$
is not differentiable.
Also it should be noted (see, \cite{Arkoun2011}) that
for the model \eqref{sec:In.1}, the adaptive pointwise estimation is possible only
in the sequential analysis framework.
That is why we study sequential estimation methods
 for the smooth
function $S$. 
In this paper we consider  the quadratic risk 
 defined as
 \begin{equation}\label{sec:In.2_risk}
\cR_\zs{p}(\wh{S}_\zs{n},S)\,=\,\E_\zs{p,S}\|\wh{S}_\zs{n}-S\|^2\,,
\quad
\|S\|^2=\int^{b}_\zs{a}\,S^2(x)\d x\,,
\end{equation}
where $\wh{S}_\zs{n}$ is an estimator of 
$S$ based on observations $(y_\zs{k})_\zs{1\le k\le n}$ and $\E_\zs{p,S}$ is the expectation with respect to the distribution 
 law $\P_\zs{p,S}$ of the process $(y_\zs{k})_\zs{1\le k\le n}$ given the distribution density $p$ and the coefficient $S$.
Moreover, taking into account that the distribution $p$ is unknown,  
we use the robust nonparametric estimation approach proposed in  \cite{GaltchoukPergamenshchikov2006a}. To this end 
 we set the robust risk  as
 \begin{equation}\label{sec:In.2_risk_rob}
\cR^{*}(\wh{S}_\zs{n},S)\,=\,
\sup_\zs{p\in\cP}
\cR_\zs{p}(\wh{S}_\zs{n},S)\,,
\end{equation}
where $\cP$ is a family of the distributions defined in Section \ref{sec:Mncnds}.

In order to estimate the function $S$
in
model
\eqref{sec:In.5} we make use of the estimator family
 $(\wh{S}_\zs{\lambda}\,,\lambda\in\Lambda)$,
where $\wh{S}_\zs{\lambda}$ is a weighted least square
estimator
with the Pinsker weights. For this family, similarly to \cite{GaltchoukPergamenshchikov2009a},
we construct
 a special 
selection rule, i.e. a
random variable $\wh{\lambda}$ with values in $\Lambda$,
 for
which we define the selection estimator as  
$\wh{S}_\zs{*}=\wh{S}_\zs{\wh{\lambda}}$. 
 Our goal in this paper is to show the non asymptotic sharp oracle inequality
for the robust risks \eqref{sec:In.2_risk_rob}, i.e. to show that  for any $\check{\varrho}>0$
and $n\ge 1$
\begin{equation}\label{sec:In.3}
\cR^{*}(\wh{S}_\zs{*},S)\,
\le\,(1+\check{\varrho})\,\min_\zs{\lambda\in\Lambda}\,
\cR^{*}(\wh{S}_\zs{\lambda},S)+
\frac{\cB_\zs{n}}{\check{\varrho} n}\,,
\end{equation}
where 
 $\cB_\zs{n}$ is a rest term such that
  for any  $\check{\delta}>0$,
\begin{equation}\label{sec:In.4}
\lim_\zs{n\to\infty}\frac{\cB_\zs{n}}{n^{\check{\delta}}}=0
\,.
\end{equation}
In this case the estimator $\wh{S}_\zs{*}$ is called optimal in the oracle inequality sense.

\noindent
In this paper, in order to obtain this inequality for model \eqref{sec:In.1}
 we develop a new model selection method
based on the truncated sequential procedures developed in
\cite{ArkounPergamenchtchikov2016}
for the pointwise efficient estimation.
Then we use the non asymptotic analysis tool proposed in \cite{GaltchoukPergamenshchikov2009a} based
on the non-asymptotic studies from
 \cite{BaronBirgeMassart1999} for a family of least-squares estimators 
and extended in \cite{FourdrinierPergamenshchikov2007}
to some other estimator families. 
To this end
we use the approach proposed in
\cite{GaltchoukPergamenshchikov20011}, i.e.
 we pass 
 to a discrete time regression model
by making use of the truncated sequential procedure introduced in 
  \cite{ArkounPergamenchtchikov2016}. 
  To this end, at any point $(z_\zs{l})_\zs{1\le l\le d}$ of a
partition of the interval $[a,b]$, we define a sequential procedure $(\tau_\zs{l},S^*_\zs{l})$ 
with a stopping rule $\tau_\zs{l}$ and an estimator $S^*_\zs{l}$. For 
$Y_\zs{l}=S^*_\zs{l}$ with $1\le l\le d$, 
we come to the regression equation on some set $\Gamma\subseteq \Omega$:
 \begin{equation}\label{sec:In.5}
Y_\zs{l}=S(z_\zs{l})+\zeta_\zs{l}\,,
\quad 1\le l\le d\,.
\end{equation}
Here,
in contrast with the classical regression model, the noise sequence 
$(\zeta_\zs{l})_\zs{1\le l\le d}$ has a complex structure, namely, 
\begin{equation}\label{sec:In.6}
\zeta_\zs{l}=\,\xi^{*}_\zs{l}+\varpi_\zs{l}\,,
\end{equation}
where $(\xi^{*}_\zs{l})_\zs{1\le l\le d}$ is a "main noise"  sequence
of  uncorrelated random variables
and
 $(\varpi_\zs{l})_\zs{1\le l\le n}$ is a sequence of
bounded random variables.

We will use the oracle inequality \eqref{sec:In.3} to 
prove the asymptotic efficiency for the proposed
procedure, using the same method as it is been used in \cite{GaltchoukPergamenshchikov2009b}.
The asymptotic efficiency means that the 
procedure provides the optimal convergence rate and 
the asymptotically minimal rate normalized risk which coincides with the Pinsker constant. It should be emphasized that
only sharp oracle inequalities of type  \eqref{sec:In.3} allow to synthesis  efficiency property in the adaptive setting.

The paper is organized as follows: In Section \ref{sec:Mncnds} we state the main conditions for the model \eqref{sec:In.1}.
In Section \ref{sec:Dt} we describe the passage to the regression scheme. In Section \ref{sec:Ms}
we describe the sequential model selection procedure. In Section \ref{sec:Mrs} 
we announce the main results.
In Section \ref{sec:Prsm} we study the properties of the obtained regression model \eqref{sec:In.5}.   
In Section \ref{sec:proofs++} we prove all basic results.  In Appendix \ref{sec:A} we give all the auxiliary and technical tools.

\medskip

\section{Main Conditions}\label{sec:Mncnds}
As in
\cite{ArkounPergamenchtchikov2016}
we assume that  in the model \eqref{sec:In.1} the i.i.d. random
variables  $(\xi_\zs{k})_\zs{k\ge 1}$  have a density  $p$
(with respect to Lebesgue measure)  from the functional class $\cP$
defined as
\begin{align}\nonumber
\cP:=&
\left\{
p\ge 0\,:\,\int^{+\infty}_\zs{-\infty}\,p(x)\,\d x=1\,,\quad
\int^{+\infty}_\zs{-\infty}\,x\,p(x)\,\d x= 0 \,,\right.
\\[2mm] \label{2.1}
&\quad\left.
\int^{+\infty}_\zs{-\infty}\,x^2\,p(x)\,\d x= 1
\quad\mbox{and}\quad
\sup_\zs{k\ge 1}\,\frac{\int^{+\infty}_\zs{-\infty}\,|x|^{2k} \,p(x)\,\d x}{\varsigma^{k}(2k-1)!!}
 \;\le 1
\,
\right\}\,,
\end{align}
where $\varsigma\ge 1$ is some fixed parameter, which may be a function
of the number observation $n$, i.e. $\varsigma=\varsigma(n)$, such that
 for any $\check{\delta}>0$
\begin{equation}
\label{varsigma-cond-1-00}
\lim_\zs{n\to\infty}
\frac{\varsigma(n)}{n^{\check{\delta}}}
=0
\,.
\end{equation}

  Note that the $(0,1)$-Gaussian density
belongs to $\cP$. In the sequel we denote this density by $ p_\zs{0}$.
It is clear that
for any $q>0$
\begin{equation}\label{sec:Sp.0}
\m^{*}_\zs{q}=\sup_\zs{p\in\cP}\,\E_\zs{p}\,|\xi_\zs{1}|^{q}
<\infty\,,
\end{equation}
where $\E_\zs{p}$ is the expectation with respect to the density $p$ from
$\cP$.
\noindent  To obtain the stable (uniformly with respect to the function $S$ ) model \eqref{sec:In.1}, we assume that for some fixed $0<\e<1$ and $L>0$ the unknown function $S$ belongs to
the $\varepsilon$ - {\em stability set} introduced in \cite{ArkounPergamenchtchikov2016}
as
\begin{equation}\label{sec:Sp.1}
\Theta_\zs{\e,L} = \left\{S\in \C_\zs{1}([a,b],\bbr) : |S|_\zs{*} \le 1-\e
\quad\mbox{and}\quad
|\dot{S}|_\zs{*}\le L
 \right\}\,,
\end{equation}
where
$\C_\zs{1}([a,b],\bbr)$ is the Banach space of continuously differentiable
$[a,b]\to\bbr$ functions and $|S|_\zs{*} = \sup_\zs{a\le x \leq b}|S(x)|$.

\section{Passage to a discrete time regression model}\label{sec:Dt}

We will use as a basic procedure the pointwise procedure from
\cite{ArkounPergamenchtchikov2016} at the points $(z_\zs{l})_\zs{1\le l\le d}$
defined as
\begin{equation}\label{sec:Dt.00-00}
z_\zs{l}=a+\frac{l}{d}(b-a)
\quad\mbox{and}\quad
d=[\sqrt{n}]
\,,
\end{equation}
where $[a]$ is the integer part of a number $a$.
 So we propose to use the first $\iota_\zs{l}$ observations for the auxiliary
estimation of $S(z_\zs{l})$. We set
\begin{equation}\label{sec:Sp.3}
\wh{S}_\zs{\iota_\zs{l}}=\frac{1}{A_\zs{\iota_\zs{l}}}\,
\sum^{\iota_\zs{l}}_\zs{j=1}\,Q_\zs{l,j}
\,y_\zs{j-1}\,y_\zs{j}\,,\quad
A_\zs{\iota_\zs{l}}=\sum^{\iota_\zs{l}}_\zs{j=1}\,Q_\zs{l,j}\,y_\zs{j-1}^2\,,
\end{equation}
where  $Q_\zs{l,j}=Q(u_\zs{l,j})$ and the kernel $Q(\cdot)$ is the indicator function of the interval $[-1;1]$, i.e. 
$Q(u)=\Chi_\zs{[-1,1]}(u)$. The points $(u_\zs{l,j})$ are defined as
\begin{equation}\label{sec:Sp.3-0-000}
u_\zs{l,j} = \frac{x_\zs{j}-z_\zs{l}}{h}
\,.
\end{equation}
Note that
 to estimate $S(z_\zs{l})$
on the basis of kernel estimator with the
kernel $Q$ we use only the observations
$(y_\zs{j})_\zs{k_\zs{1,l}\le j\le k_\zs{2,l}}$
from the $h$ -
neighborhood
of the point $z_\zs{l}$, i.e.
\begin{equation}\label{sec:Sp.7}
k_\zs{1,l}= [n\wt{z}_\zs{l} - n\wt{h}] + 1
\quad\mbox{and}\quad
k_\zs{2,l}=  [n\wt{z}_\zs{l} + n\wt{h}]
\wedge n
\,,
\end{equation}
where $\wt{z}_\zs{l}=(z_\zs{l}-a)/(b-a)$ and $\wt{h}=h/(b-a)$.
Note that, only for the last point $z_\zs{d}=b$, we take $k_\zs{2,d}=n$.
We chose $\iota_\zs{l}$
in \eqref{sec:Sp.3}
 as
\begin{equation}\label{sec:Sp.9-0}
\iota_\zs{l}=k_\zs{1,l}+\q
\quad\mbox{and}\quad
\q=\q_\zs{n}=[(n\wt{h})^{\mu_\zs{0}}]
\end{equation}
for some $0<\mu_\zs{0}<1$.
In the sequel for any $0\le k< m\le n$ we set
\begin{equation}\label{sec:Sp.3-0}
A_\zs{k,m}=\sum^{m}_\zs{j=k+1}\,Q_\zs{l,j}\,y_\zs{j-1}^2
\quad\mbox{and}\quad
A_\zs{m}=A_\zs{0,m}
\,.
\end{equation}
\noindent 
Next, similarly to \cite{Arkoun2011},
we use a kernel sequential procedure
 based on the observations
$(y_\zs{j})_\zs{\iota_\zs{l}\le j\le n}$.
To transform the kernel estimator in a linear function of observations and
we replace the number of observations $n$ by the following stopping time
\begin{equation}\label{sec:Sp.4}
 \tau_\zs{l} = \inf\{\iota_\zs{l} +1\le k\le k_\zs{2,l}\,: \,A_\zs{\iota_\zs{l},k} \ge H_\zs{l} \}
 \,,
 \end{equation}
where $\inf\{\emptyset\}=k_\zs{2,l}$ and the positive threshold $H_\zs{l}$ will be chosen as a positive
random variable which is measurable with respect to the $\sigma$-field
$\{y_\zs{1},\ldots,y_\zs{\iota_\zs{l}}\}$.  

\noindent
Now we define the sequential estimator as

\begin{equation}\label{sec:Sp.5}
S^{*}_\zs{l}=\frac{1}{H_\zs{l}}\,
\left(\sum^{\tau_\zs{l}-1}_\zs{j=\iota_\zs{l}+1}\,Q_\zs{l,j}\,y_\zs{j-1}\,y_\zs{j}\,+\,\varkappa_\zs{l}\,Q_\zs{l,\tau_\zs{l}}\,y_\zs{\tau_\zs{l}-1}\,y_\zs{\tau_\zs{l}}\right)
\Chi_\zs{\Gamma_\zs{l}}\,,
\end{equation}
where
$\Gamma_\zs{l}=\{A_\zs{\iota_\zs{l},k_\zs{2,l}-1}\ge  H_\zs{l}\}$ and
 the correcting  coefficient $0<\varkappa_\zs{l}\le 1$ on this set
is defined as
\begin{equation}\label{sec:Sp.5-cr}
 A_\zs{\iota_\zs{l},\tau_\zs{l}-1} + \varkappa^{2}_\zs{l}\,Q_\zs{l,\tau_\zs{l}}
y^2_\zs{\tau_\zs{l}-1} = H_\zs{l}
\,.
\end{equation}

\noindent
Note that, to obtain the efficient kernel estimator of $S(z_\zs{l})$ we need to use all $k_\zs{2,l}-\iota_\zs{l}-1$ observations.
Similarly to  \cite{KonevPergamenshchikov1984}, one can show that $\tau_\zs{l}\approx \gamma_\zs{l}\,H_\zs{l}$ as $H_\zs{l}\to\infty$, where
\begin{equation}\label{sec:Sp.6}
\gamma_\zs{l}=1-S^2(z_\zs{l})\,.
\end{equation}
\noindent 
Therefore, one needs to chose $H_\zs{l}$ as $(k_\zs{2,l}-\iota_\zs{l}-1)/\gamma_\zs{l}$.
Taking into account  that the coefficients $\gamma_\zs{l}$ are unknown,
we define the threshold $H_\zs{l}$ as

\begin{equation}\label{sec:Sp.8}
H_\zs{l}=\frac{1-\wt{\epsilon}}{\wt{\gamma}_\zs{l}}
(k_\zs{2,l}-\iota_\zs{l}-1)
\quad\mbox{and}\quad 
\wt{\epsilon}=\frac{1}{2+\ln n}
\,,
\end{equation}
where $\wt{\gamma}_\zs{l}=1-\wt{S}^{2}_\zs{\iota_\zs{l}}$ and $\wt{S}_\zs{\iota_\zs{l}}$ is the projection of the estimator
$\wh{S}_\zs{\iota_\zs{l}}$ in the interval $]-1+\wt{\e},1-\wt{\e}[$, i.e.
\begin{equation}
\label{est-wt-S-00}
\wt{S}_\zs{\iota_\zs{l}} = \min (\max (\wh{S}_\zs{\iota_\zs{l}},-1+\wt{\e}),1-\wt{\e})\,.
\end{equation}
To obtain the uncorrelated stochastic terms 
in kernel estimator for $S(z_\zs{l})$
 we chose the bandwidth $h$ as
\begin{equation}\label{sec:Sp.9}
h=\frac{b-a}{2 d}
\,.
\end{equation}
As to the estimator $\wh{S}_\zs{\iota_\zs{l}}$, we can show the following property.

\begin{proposition}\label{Pr.sec:Prs.stp.times.1}
The convergence rate
in probability 
 of the estimator \eqref{est-wt-S-00}
is more rapid than any power function, i.e.
for any $\b>0$ 
\begin{equation}\label{sec:Prs.stp.times.1}
\lim_\zs{n\to\infty}\,n^{\b}
\max_\zs{1\le l\le d}
\sup_\zs{S\in \Theta_\zs{\e,L}}
\sup_\zs{p\in \cP}
\P_\zs{p,S}\left(\vert \wt{S}_\zs{\iota_\zs{l}}-S(z_\zs{l})\vert>\epsilon_\zs{0}\right)
=0\,,
\end{equation}
where $\epsilon_\zs{0}=\epsilon_\zs{0}(n)\to 0$ as $n\to\infty$ such that 
$\lim_\zs{n\to\infty}n^{\check{\delta}}\epsilon_\zs{0}=\infty$
for any $\check{\delta}>0$.
\end{proposition}

\noindent 
Now we set
\begin{equation}\label{sec:In.df_est}
Y_\zs{l}=S_\zs{l}^{*}(z_\zs{l})\Chi_\zs{\Gamma}
\quad\mbox{and}\quad
\Gamma=
\cap^{d}_\zs{l=1}
\,\Gamma_\zs{l} 
\,.
\end{equation}

\noindent Using the convergence \eqref{sec:Prs.stp.times.1},
we  study the  probability properties of the set $\Gamma$
 in the  following proposition.

\begin{proposition}\label{Pr.sec:Prs.stp.times.22}
For any $\b>0$, the probability of the set $\Gamma$ satisfies the following
asymptotic equality
\begin{equation}\label{sec:gamma_lim}
\lim_\zs{n\to\infty}\,n^{\b}
\sup_\zs{S\in \Theta_\zs{\e,L}}
\P_\zs{p,S}\left(\Gamma^{c}\right)
=0\,.
\end{equation}
\end{proposition}

\noindent
In view of this proposition  we can shrink the set $\Gamma^{c}$. So,
using the estimators \eqref{sec:In.df_est}
on the set $\Gamma$ we obtain the discrete time regression model
\eqref{sec:In.5}
in which
\begin{equation}
\label{sec:In.5++}
\xi^{*}_\zs{l}
=\frac{
\sum^{\tau_\zs{l}-1}_\zs{j=\iota_\zs{l}+1}\,Q_\zs{l,j}\,y_\zs{j-1}\,\xi_\zs{j}+\,\varkappa_\zs{l}\,Q(u_\zs{l,\tau_\zs{l}})\,y_\zs{\tau_\zs{l}-1}\,\xi_\zs{\tau_\zs{l}}
}{H_\zs{l}}
\end{equation}
and
$\varpi_\zs{l}=\varpi_\zs{1,l}+\varpi_\zs{2,l}$,
where
$$
\varpi_\zs{1,l}=
\frac{
\sum^{\tau_\zs{l}-1}_\zs{j=\iota_\zs{l}+1}\,Q_\zs{l,j}\,y^{2}_\zs{j-1}\check{s}_\zs{l,j}+
\varkappa^{2}_\zs{l}\,Q(u_\zs{l,\tau_\zs{l}})\,y^{2}_\zs{\tau_\zs{l}-1}
\check{s}_\zs{l,\tau_\zs{l}}
}{H_\zs{l}}\,,\quad
\check{s}_\zs{l,j}=S(x_\zs{j})-S(z_\zs{l})
$$
and
$$
\varpi_\zs{2,l}=
\frac{
(\varkappa_\zs{l}-\varkappa^{2}_\zs{l})\,Q(u_\zs{l,\tau_\zs{l}})\,y^{2}_\zs{\tau_\zs{l}-1}\,S(x_\zs{\tau_\zs{l}})
}{H_\zs{l}}
\,.
$$
Note that in the model
\eqref{sec:In.5}
the random variables 
$(\xi^{*}_\zs{j})_\zs{1\le j\le d}$
are defined only on the set $\Gamma$. For technical reasons we need to define these variables on the
set $\Gamma^{c}$ as well.  
To this end, for any $j\ge1$ we set
\begin{equation}\label{sec:Sp.3-0++}
\check{Q}_\zs{l,j}=Q_\zs{l,j}\,y_\zs{j-1}\,\Chi_\zs{\{j< k_\zs{2,l}\}}
+\sqrt{H_\zs{l}}\,Q_\zs{l,j}\,\Chi_\zs{\{j=k_\zs{2,l}\}}
\end{equation}
and
$\check{A}_\zs{\iota_\zs{l},m}=\sum^{m}_\zs{j=\iota_\zs{l}+1}\,\check{Q}^{2}_\zs{l,j}$. Note, that
for any $j\ge 1$ and $l\neq m$
\begin{equation}\label{sec:def_eta++}
\check{Q}_\zs{l,j}\,\check{Q}_\zs{m,j}\,
=0
\,.
\end{equation}
and
 $\check{A}_\zs{\iota_\zs{l},k_\zs{2,l}}\ge H_\zs{l}$. So now we can modify the stopping time
\eqref{sec:Sp.4} as

\begin{equation}\label{sec:Sp.4+++}
 \check{\tau}_\zs{l} = \inf\{ k\geq \iota_\zs{l}+1: \check{A}_\zs{\iota_\zs{l},k} \ge H_\zs{l} \}\,.
 \end{equation}
Obviously,  $\check{\tau}_\zs{l}\le k_\zs{2,l}$  and 
$\check{\tau}_\zs{l}=\tau_\zs{l}$ on the set $\Gamma$ for any $1\le l\le d$.
Now similarly to
\eqref{sec:Sp.5-cr} we define the correction coefficient as
\begin{equation}\label{sec:Sp.5-cr++}
 \check{A}_\zs{\iota_\zs{l},\check{\tau}_\zs{l}-1} + \check{\varkappa}^{2}_\zs{l}\,
 \check{Q}^{2}_\zs{l,\check{\tau}_\zs{l}}
 = H_\zs{l}
\,.
\end{equation}
It is clear that $0< \check{\varkappa}_\zs{l} \le 1$
and $\check{\varkappa}_\zs{l}=\varkappa_\zs{l}$
on the set $\Gamma$ for $1\le l\le d$. Using this coefficient we set
\begin{equation}\label{sec:def_eta}
\eta_\zs{l}
=\frac{
\sum^{\check{\tau}_\zs{l}-1}_\zs{j=\iota_\zs{l}+1}\,\check{Q}_\zs{l,j}\,
\xi_\zs{j}+\,\check{\varkappa}_\zs{l}\,\check{Q}_\zs{l,\check{\tau}_\zs{l}}\,\xi_\zs{\check{\tau}_\zs{l}}}{H_\zs{l}}
\,.
\end{equation}
Note that on the set $\Gamma$, for any $1\le l\le d$, the random variables $\eta_\zs{l}=\xi^{*}_\zs{l}$.
Moreover (see Lemma \ref{Le.sec:App.eta_prps}), 
 for any $1\le l\le d$ and $p\in\cP$
\begin{equation}\label{sec:def_eta++_prs_0}
\E_\zs{p,S}\left(\eta_\zs{l}\,\vert \cG_\zs{l}\right)=0
\,,
\quad
\E_\zs{p,S}\left(\eta^{2}_\zs{l}\,\vert \cG_\zs{l}\right)=\sigma^{2}_\zs{l}
\quad\mbox{and}\quad
\E_\zs{p,S}\left(\eta^{4}_\zs{l}\,\vert \cG_\zs{l}\right)\,\le \check{\m}\sigma^{4}_\zs{l}
\,,
\end{equation}
where  $\sigma_\zs{l}=H^{-1/2}_\zs{l}$, $\cG_\zs{l}=\sigma\{\eta_\zs{1},\ldots,\eta_\zs{l-1},\sigma_\zs{l}\}$ 
 and $\check{\m}=4(144/\sqrt{3})^{4}\,\m^{*}_\zs{4}$.
It is clear that
\begin{equation}
\label{sec:seq_bound_sigma}
\sigma_\zs{0,*}
\le
\min_\zs{1\le l\le d}\sigma^{2}_\zs{l}
\le 
\max_\zs{1\le l\le d}\sigma^{2}_\zs{l}
\le \sigma_\zs{1,*}
\,,
\end{equation}
where
$$
\sigma_\zs{0,*}=
\frac{1-\epsilon^{2}}{2(1-\wt{\epsilon})nh}
\quad\mbox{and}\quad
\sigma_\zs{1,*}=
\frac{1}{(1-\wt{\epsilon})(2nh-\q-3)}
\,.
$$

\noindent 
Now, taking into account that
$
\vert \varpi_\zs{1,l}\vert\le L h$,
 for any $S\in \Theta_\zs{\e,L} $
we obtain that
\begin{equation}\label{sec:bound_varpi}
\sup_\zs{S\in \Theta_\zs{\e,L}}\,
\E_\zs{p,S} \Chi_\zs{\Gamma} \varpi^{2}_\zs{l}
\le 
 \left(
 L^{2} h^{2}+\frac{\check{\upsilon}_\zs{n}}{(n h)^{2}}
 \right)
 \,, 
\end{equation}
where $\check{\upsilon}_\zs{n}=\sup_\zs{p\in \cP}\,\sup_\zs{S\in \Theta_\zs{\e,L}}\,\E_\zs{p,S}\,\max_\zs{1\le j\le n}\, y^{4}_\zs{j}$.
The behavior of this coefficient is studied in the following proposition.

\begin{proposition}\label{Pr.sec:Ms.0}
For any $\b>0$ the sequence $(\check{\upsilon}_\zs{n})_\zs{n\ge 1}$
satisfies the following limiting equality
\begin{equation}\label{sec:yyy_lim}
\lim_\zs{n\to\infty}\,n^{-\b}
\,
 \check{\upsilon}_\zs{n}
=0\,.
\end{equation}
\end{proposition}

\begin{remark}
\label{Re.sec:seq-pr}
It should be noted that the property 
\eqref{sec:yyy_lim} means that the asymptotic behavior of the upper bound \eqref{sec:bound_varpi} is
approximately as $h^{-2}$ when $n\to\infty$. We will use this in the oracle inequalities below.
\end{remark}

\begin{remark}
\label{Re.sec:seq-pr-0-1}
It should be emphasized that
to estimate the function $S$ in 
\eqref{sec:In.1}
 we use the approach developed in
\cite{GaltchoukPergamenshchikov20011}
for the sequential drift estimation problem
in the stochastic differential equation. On the basis
of the efficient sequential kernel
procedure developped in
\cite{GaltchoukPergamenshchikov2005}, \cite{GaltchoukPergamenshchikov2006b} and \cite{GaltchoukPergamenshchikov2015}
with the kernel-indicator, the stochastic differential equation 
is replaced by regression model. It should be noted that to obtain the efficient
estimator one needs to take the kernel-indicator estimator.
By this reason, in this paper, we use 
the kernel-indicator in the sequential estimator
\eqref{sec:Sp.5}. It also should be noted that the 
sequential estimator \eqref{sec:Sp.5} which has the  same
form as in 
\cite{ArkounPergamenchtchikov2016}, except the last term, in which the correction coefficient is replaced
by the square root of the coefficient used in  \cite{Konev2016}.
We modify this procedure
to calculate the variance of the stochastic term  \eqref{sec:In.5++}. 
\end{remark}

\section{Model selection}\label{sec:Ms}

In this section we consider the nonparametric estimation problem
in the non asymptotic setting
 for the 
regression model \eqref{sec:In.5} for some set $\Gamma\subseteq\Omega$.
 The design points
$(z_\zs{l})_\zs{1\le l\le d}$
are defined in \eqref{sec:Dt.00-00}. The function $S(\cdot)$ is unknown and has to be estimated
from observations the $Y_1,\ldots,Y_\zs{d}$. Moreover, we assume that  the unobserved random variables $(\eta_\zs{l})_\zs{1\le l\le d}$ satisfy the properties
\eqref{sec:def_eta++_prs_0} with some  nonrandom constant $\check{\m}>1$ and the known random positive coefficients
  $(\sigma_\zs{l})_\zs{1\le l\le d}$ satisfy the inequality \eqref{sec:seq_bound_sigma}
  for some nonrandom positive constants $\sigma_\zs{0,*}$ and $\sigma_\zs{1,*}$
Concerning the random sequence $\varpi=(\varpi_\zs{l})_\zs{1\le l\le n}$ we suppose that 
\begin{equation}\label{sec:Ms.3-1}
\u^{*}_\zs{d}=
\E_\zs{p,S}\Chi_\zs{\Gamma}\|\varpi\|^2_\zs{d}<\infty\,.
\end{equation}
The performance of any estimator $\wh{S}$ will be measured by
the empirical squared error
$$
\|\wh{S}-S\|^2_\zs{d}=(\wh{S}-S,\wh{S}-S)_\zs{d}
=\frac{b-a}{d}\sum^{d}_\zs{l=1}(\wh{S}(z_\zs{l})-S(z_\zs{l}))^2\,.
$$
Now we fix a basis 
$(\phi_\zs{j})_\zs{1\le j\le n}$ which is orthonormal for the empirical inner product:
\begin{equation}\label{sec:Ms.4nn_1}
(\phi_i\,,\,\phi_\zs{j})_\zs{d}=
\frac{b-a}{d}\sum^{d}_\zs{l=1}\,\phi_i(z_\zs{l})\phi_\zs{j}(z_\zs{l})= \Chi_\zs{\{i=j\}}
\,.
\end{equation}
For example, we can take the trigonometric basis  $(\phi_\zs{j})_\zs{j\ge\, 1}$ in $\L_\zs{2}[a,b]$ defined as
\begin{equation}\label{sec:In.5_trg}
\phi_\zs{1} = 1\,,\quad \phi_\zs{j}(x)= \sqrt{\frac{2}{b-a}} \Trg_\zs{j}\left(2\pi[j/2]\l_\zs{0}(x)\right)\,,\quad j\ge\,2\,,
\end{equation}
where the function $\Trg_\zs{j}(x)= \cos(x)$ for even $j$ and $\Trg_\zs{j}(x)= \sin(x)$ for odd $j$, $[x]$ denotes the integer part of $x$. 
and $\l_\zs{0}(x)=(x-a)/(b-a)$.
Note that, using the orthonormality property \eqref{sec:Ms.4nn_1}
 we can represent for any $1\le l\le d$ the function $S$ as
\begin{equation}\label{sec:In.5_s--1}
S(z_\zs{l})=\sum^{d}_\zs{j=1}\,\theta_\zs{j,d}\,\phi_\zs{j}(z_\zs{l})
\quad\mbox{and}\quad
\theta_\zs{j,d}=
\left(S,\phi_\zs{j} \right)_\zs{d}
\,.
\end{equation}
So, to estimate the function $S$ we have to estimate the Fourier coefficients $(\theta_\zs{j,d})_\zs{1\le j\le d}$. To this end, 
we replace the function $S$ by these observations, i.e.
\begin{equation}\label{sec:In.5_est--1}
\wh{\theta}_\zs{j,d}=\frac{b-a}{d}
\sum^{d}_\zs{l=1}\,Y_\zs{l}\phi_\zs{j}(z_\zs{l})\,.
\end{equation}
From \eqref{sec:In.5} we obtain immediately the following regression scheme 
\begin{equation}\label{sec:Ms.5}
\wh{\theta}_\zs{j,d}\,=\,\theta_\zs{j,d}\,+\,\zeta_\zs{j,d}
\quad\mbox{with}\quad
\zeta_\zs{j,d}=\sqrt{\frac{b-a}{d}}
\eta_\zs{j,d}+\varpi_\zs{j,d}\,,
\end{equation}
where
$$
\eta_\zs{j,d}=\sqrt{\frac{b-a}{d}}\sum^{d}_\zs{l=1}\,
\eta_\zs{l}\phi_\zs{j}(z_\zs{l})
\quad\mbox{and}\quad
\varpi_\zs{j,d}=\frac{b-a}{d}
\sum^{d}_\zs{l=1}\,\varpi_\zs{l}\,\phi_\zs{j}(z_\zs{l})\,.
$$
Note that the upper bound \eqref{sec:seq_bound_sigma} and the Bounyakovskii-Cauchy-Schwarz inequality imply that
\begin{equation}\label{sec:Ms.6}
|\varpi_\zs{j,d}|\le\|\varpi\|_\zs{d}\,\|\phi_\zs{j}\|_\zs{d}=
\|\varpi\|_\zs{d}\,.
\end{equation}
We estimate the function $S$ on the grid 
\eqref{sec:Dt.00-00}
by the weighted least-squares estimator
\begin{equation}\label{sec:Ms.7}
\wh{S}_\zs{\lambda}(z_\zs{l})\,=\,\sum^{d}_\zs{j=1}\,\lambda(j)\,\wh{\theta}_\zs{j,d}\,\phi_\zs{j}(z_\zs{l})\,
\Chi_\zs{\Gamma}\,,\quad 1\le l\le d\,,
\end{equation}
where the weight vector $\lambda=(\lambda(1),\ldots,\lambda(d))'$ 
belongs to some finite set $\Lambda\subset [0,1]^{d}$, the prime denotes the transposition.
We set for any $a\le t\le b$
\begin{equation}\label{sec:Ms.8}
\wh{S}_\zs{\lambda}(t)=
\wh{S}_\zs{\lambda}(z_\zs{1})
\Chi_\zs{\{a\le t\le z_\zs{1}\}}+
\sum^{d}_\zs{l=2}\wh{S}_\zs{\lambda}(z_\zs{l})
\Chi_\zs{\{z_\zs{l-1}< t\le z_\zs{l}\}}\,.
\end{equation}
\noindent Moreover, denoting 
$\lambda^2=(\lambda^2(1),\dots,\lambda^2(n))'$
we define the following sets
\begin{equation}\label{sec:Ms.9}
\Lambda_\zs{1}=\{\lambda^2\,,\,\lambda\in\Lambda\}
\quad\mbox{and}\quad
\Lambda_\zs{2}=\Lambda\cup\Lambda_\zs{1}\,.
\end{equation}

\noindent Denote by $\nu$ the cardinal number of the set $\Lambda$ and
$$
\nu^*=\max_\zs{\lambda\in\Lambda}\,
\sum^{d}_\zs{j=1}
\Chi_\zs{\{\lambda(j)>0\}}\,.
$$
\noindent In order to obtain
a good estimator, we have to write a rule to choose a weight vector 
$\lambda\in\Lambda$ in \eqref{sec:Ms.7}.  We define
 the empirical squared risk as
$$
\Er_\zs{d}(\lambda)=\|\wh{S}_\zs{\lambda}-S\|^2_\zs{d}\,.
$$
Using \eqref{sec:In.5_s--1} and \eqref{sec:Ms.7} we can rewrite this risk as
\begin{equation}\label{sec:Ms.10}
\Er_\zs{d}(\lambda)=
\sum^{d}_\zs{j=1}\,\lambda^2(j)\wh{\theta}^2_\zs{j,d}\,-
2\,\sum^{d}_\zs{j=1}\,\lambda(j)\wh{\theta}_\zs{j,d}\,\theta_\zs{j,d}\,+\,
\sum^{d}_\zs{j=1}\,\theta^2_\zs{j,d}\,.
\end{equation}
Since the coefficient $\theta_\zs{j,d}$ is unknown, we need to replace the term 
$\wh{\theta}_\zs{j,d}\,\theta_\zs{j,d}$
by some of its estimators which we choose as 
\begin{equation}\label{sec:Ms.11}
\wt{\theta}_\zs{j,d}=
\wh{\theta}^2_\zs{j,d}-\frac{b-a}{d}\,\s_\zs{j,d}
\quad\mbox{with}\quad
\s_\zs{j,d}=\frac{b-a}{d}\,\sum^{d}_\zs{l=1}\,\sigma^2_\zs{l}\,\phi^2_\zs{j}(z_\zs{l})\,.
\end{equation}
Note that from 
\eqref{sec:seq_bound_sigma}
-
\eqref{sec:Ms.4nn_1}
 it follows that
\begin{equation}\label{sec:Ms.11++0}
\s_\zs{j,d}\le \sigma_\zs{1,*}\,.
\end{equation}
Finally, we define the cost function of the form
\begin{equation}\label{sec:Ms.12}
J_\zs{d}(\lambda)\,=\,\sum^{d}_\zs{j=1}\,\lambda^2(j)\wh{\theta}^2_\zs{j,d}\,-
2\,\sum^{d}_\zs{j=1}\,\lambda(j)\,\wt{\theta}_\zs{j,d}\,
+\,\delta P_\zs{d}(\lambda)\,,
\end{equation}
where  the penalty term is defined as
\begin{equation}\label{sec:pen_term}
P_\zs{d}(\lambda)=
\frac{b-a}{d}
\sum^{d}_\zs{j=1}\lambda^2(j)\s_\zs{j,d}
\end{equation}
and $0< \delta< 1$ is some positive constant which will be chosen later.
We set
\begin{equation}\label{sec:Ms.13}
\wh{\lambda}=\mbox{argmin}_\zs{\lambda\in\Lambda}\,J_\zs{d}(\lambda)
\end{equation}
and  define an estimator of $S(t)$ of the form \eqref{sec:Ms.8}:
\begin{equation}\label{sec:Ms.14}
\wh{S}_\zs{*}(t)=\wh{S}_\zs{\wh{\lambda}}(t)
\quad\mbox{for}\quad a\le t\le b\,.
\end{equation}

\begin{remark}
\label{Re.sec:ModSelec-00+}
We 
use the procedure \eqref{sec:Ms.14}
to estimate the function $S$ in the autoregressive model
\eqref{sec:In.1}
 through the regression scheme  \eqref{sec:In.5}
generated by the sequential procedures 
\eqref{sec:In.df_est}.
\end{remark}

\section{Main results}\label{sec:Mrs}

In this section we formulate all main results. First we obtain the sharp oracle inequality for the selection model procedure
\eqref{sec:Ms.14} for the general regression model \eqref{sec:In.5}.

\begin{theorem}\label{Th.sec:Ms.1-0+++}
There exists some constant $\l^{*}>0$ such that for any 
weight vectors set $\Lambda$, any $p\in\cP$,
 any $n\ge\,1 $ and $ 0<\delta \leq 1/12$, 
 the procedure \eqref{sec:Ms.14}, 
satisfies the following oracle inequality
\begin{align}\nonumber
\E_\zs{p,S}\Vert\wh{S}_\zs{*}-S\Vert^{2}_\zs{d}&\leq\frac{1+4\delta}{1-6\delta} \min_\zs{\lambda\in\Lambda} 
\,
\E_\zs{p,S}\Vert\wh{S}_\zs{\lambda}-S\Vert^{2}_\zs{d}
\\[2mm] \label{sec:Mrs.1--00-22}
&+
\l^{*}
\frac{\nu\varsigma^{2}}{\delta}
\left(
\frac{\sigma^{2}_\zs{1,*}}{\sigma_\zs{0,*}d}
+\u^{*}_\zs{d}
+
\delta^2\sqrt{\P_\zs{S}(\Gamma^c)}
\right)
\,.
\end{align}
\end{theorem}

\noindent
Using now Lemma
\ref{Le.sec:A.1-1}
we obtain the oracle inequality for the quadratic risks \eqref{sec:In.2_risk}.

\begin{theorem}\label{Pr.sec:Ms.1++}
There exists some constant $\l^{*}>0$ such that for any 
weight vectors set $\Lambda$, any continuously differentiable function $S$,
any $p\in\cP$,
any $n\ge\,1 $ and $ 0 <\delta \leq 1/12$, 
 the procedure \eqref{sec:Ms.14} 
satisfies the following oracle inequality
\begin{align}\nonumber
\cR_\zs{p}(\wh{S}_\zs{*},S)
\,
&\leq\frac{(1+4\delta)(1+\delta)^2}{1-6\delta} \min_\zs{\lambda\in\Lambda} 
\,
\cR_\zs{p}(\wh{S}_\zs{\lambda},S)
\\[2mm] \label{sec:Mrs.1--00-33}
&+
\l^{*}
\frac{\varsigma^{2} \nu}{\delta}
\left(
\frac{\Vert\dot{S}\Vert^{2}}{d^{^2}}
+
\frac{\sigma^{2}_\zs{1,*}}{\sigma_\zs{0,*}d}
+\u^{*}_\zs{d}
+
\delta^2\sqrt{\P_\zs{S}(\Gamma^c)}
\right)
\,.
\end{align}
\end{theorem}

Now we assume that the cardinal $\nu$ of $\Lambda$ and the parameter $\varsigma$
in the density family
\eqref{2.1} are 
functions of the number observations $n$, i.e.
$\nu=\nu(n)$
and 
$\varsigma=\varsigma(n)$
such that for any $\check{\delta}>0$
\begin{equation}
\label{cond-1-00}
\lim_\zs{n\to\infty}
\frac{\nu(n)}{n^{\check{\delta}}}
=0
\,.
\end{equation}

\medskip
\noindent 
Using Propositions \ref{Pr.sec:Prs.stp.times.22} -- \ref{Pr.sec:Ms.0}
and the bounds
 \eqref{sec:seq_bound_sigma}
-
\eqref{sec:bound_varpi}
we obtain the oracle inequality for the estimation problem for the model \eqref{sec:In.1}.
\begin{theorem}\label{Th.sec:Ms.2--0}
Assume that the conditions \eqref{varsigma-cond-1-00}
and
\eqref{cond-1-00}
hold.
 Then for any $p\in\cP$,
$S\in\Theta_\zs{\e,L}$,
 $n\ge\, 3$ and $ 0 <\delta \leq 1/12$, 
 the procedure \eqref{sec:Ms.14} 
satisfies the following oracle inequality
\begin{equation}\label{sec:Mrs.1--00}
\cR_\zs{p}(\wh{S}_*,S)\leq\frac{(1+4\delta)(1+\delta)^2}{1-6\delta} \min_\zs{\lambda\in\Lambda} 
\cR_\zs{p}(\wh{S}_\lambda,S)
+
\,
\frac{\check{\B}_\zs{n}(p)}{\delta n}
\,,
\end{equation}
where the term $\check{\B}_\zs{n}(p)$ is such that for any $\check{\delta}>0$
$$
\lim_\zs{n\to\infty}\,\frac{\check{\B}_\zs{n}(p)}{n^{\check{\delta}}}=0\,.
$$
\end{theorem}

\noindent 
We obtain the same inequality for the robust risks

\begin{theorem}\label{Th.sec:Ms.2}
Assume that the conditions \eqref{varsigma-cond-1-00}
and
\eqref{cond-1-00}
hold.
 Then for any
 $n\ge\, 3$, any
$S\in\Theta_\zs{\e,L}$ and any $ 0 <\delta \leq 1/12$, 
 the procedure \eqref{sec:Ms.14} 
satisfies the following oracle inequality
\begin{equation}\label{sec:Mrs.1--00-12}
\cR^{*}(\wh{S}_*,S)\leq\frac{(1+4\delta)(1+\delta)^2}{1-6\delta} \min_\zs{\lambda\in\Lambda} 
\cR^{*}(\wh{S}_\lambda,S)
+
\,
\frac{\B^{*}_\zs{n}}{\delta n}
\,,
\end{equation}
where the term $\check{\B}_\zs{n}$ is such that for any $\check{\delta}>0$
$$
\lim_\zs{n\to\infty}\,\frac{\B^{*}_\zs{n}}{n^{\check{\delta}}}=0\,.
$$
\end{theorem}

It is well known that to obtain the efficiency property we need to specify the weight coefficients 
$(\lambda(j))_\zs{1\le j\le n}$ (see, for example, \cite{GaltchoukPergamenshchikov2009b}).
 Consider for some fixed $0<\varepsilon<1$
a numerical grid of the form
\begin{equation}\label{sec:Ga.0}
\cA=\{1,\ldots,k^*\}\times\{\varepsilon,\ldots,m\varepsilon\}\,,
\end{equation}
where $m=[1/\ve^2]$. We assume that both
parameters $k^*\ge 1$ and $\varepsilon$ are functions of $n$, i.e.
$k^*=k^*(n)$ and $\ve=\ve(n)$, such that
\begin{equation}\label{sec:Ga.1}
\left\{
\begin{array}{ll}
&\lim_\zs{n\to\infty}\,k^*(n)=+\infty\,,
\quad\lim_\zs{n\to\infty}\,\dfrac{k^*(n)}{\ln n}=0\,,\\[6mm]
&
\lim_\zs{n\to\infty}\,\varepsilon(n)=0
\quad\mbox{and}\quad
\lim_\zs{n\to\infty}\,n^{\check{\delta}}\ve(n)\,=+\infty
\end{array}
\right.
\end{equation}
for any $\check{\delta}>0$. One can take, for example, for $n\ge 2$
\begin{equation}\label{sec:Ga.1-00}
\ve(n)=\frac{1}{ \ln n }
\quad\mbox{and}\quad
k^*(n)=k^{*}_\zs{0}+\sqrt{\ln n}\,,
\end{equation}
where $k^{*}_\zs{0}\ge 0$ is some fixed constant.
 For each $\alpha=(\beta, \l)\in\cA$, we introduce the weight
sequence
$$
\lambda_\zs{\alpha}=(\lambda_\zs{\alpha}(j))_\zs{1\le j\le n}
$$
with the elements
\begin{equation}\label{sec:Ga.2}
\lambda_\zs{\alpha}(j)=\Chi_\zs{\{1\le j<j_\zs{*}\}}+
\left(1-(j/\omega_\alpha)^\beta\right)\,
\Chi_\zs{\{ j_\zs{*}\le j\le \omega_\zs{\alpha}\}},
\end{equation}
where
$j_\zs{*}=1+\left[\ln n\right]$, $\omega_\zs{\alpha}=(\d_\zs{\beta}\,\l\, n)^{1/(2\beta+1)}$
and
$$
\d_\zs{\beta}=\frac{(\beta+1)(2\beta+1)}{\pi^{2\beta}\beta}
\,.
$$
Now we define the set $\Lambda$ 
as
\begin{equation}\label{sec:Ga.3}
\Lambda\,=\,\{\lambda_\zs{\alpha}\,,\,\alpha\in\cA\}\,.
\end{equation}

\noindent 
Note that these weight coefficients are
  used in \cite{KonevPergamenshchikov2012, KonevPergamenshchikov2015}
 for continuous time regression
models  to show the asymptotic efficiency.
It will be noted that in this case the cardinal of the set $\Lambda$ is  
$\nu=k^{*} m$. It is clear that properties \eqref{sec:Ga.1}
imply condition \eqref{cond-1-00}.

\section{Properties of the regression model \eqref{sec:In.5}}\label{sec:Prsm}


In order to prove the oracle inequalities we  need to study  the conditions
introduced in \cite{KonevPergamenshchikov2012} for the general semi-martingale model. 
To this end we set for any $\lambda\in\bbr^{d}$ the functions
\begin{equation}\label{sec:Prsm.1}
 \B(\lambda)
=
\frac{b-a}{\sqrt{d}}\,
 \sum_\zs{j=1}^{d}\,\lambda(j)\,\wt{\eta}_\zs{j,d}
\,,
\quad
\wt{\eta}_\zs{j,d}=\eta^2_\zs{j,d}- \E_\zs{p,S}\,\eta^2_\zs{j,d}
\,.
\end{equation}

\begin{proposition}\label{Pr.sec:A.06-3}
For any $d\ge 1$ and any $\lambda=(\lambda_\zs{1},\ldots,\lambda_\zs{d})\in\bbr^{d}$
\begin{equation}\label{sec:L_2-Upp_2}
  \E_\zs{p,S}\,\Chi_\zs{\Gamma}\,
  \,
  \B^{2}(\lambda)
\,
\le 
10(b-a)
\sigma_\zs{1,*}
\check{\m}\,
\E_\zs{p,S}\,
 P_\zs{d}(\lambda)
\,.
  \end{equation}
where $\check{\m}$  
  is defined in
  \eqref{sec:def_eta++_prs_0}.
\end{proposition}

\proof
First note that the random variable $\wt{\eta}_\zs{j,d}$ can be represented as
$$
\wt{\eta}_\zs{j,d}=
\frac{b-a}{d}
\sum^{d}_\zs{l=1}
\left(
\phi^2_\zs{j}(z_\zs{l})\,\check{\eta}_\zs{l}
+
2\Chi_\zs{\{l\ge 2\}}\,
\check{\upsilon}_\zs{j,l}\eta_\zs{l}
\right)\,,
$$
where 
$\check{\eta}_\zs{l}=\eta^{2}_\zs{l}-\sigma^{2}_\zs{l}$ 
 and 
$\check{\upsilon}_\zs{j,l}=\phi_\zs{j}(z_\zs{l})
\sum^{l-1}_\zs{r=1}\,\phi_\zs{j}(z_\zs{r})\eta_\zs{r}$.
Therefore, we can rewrite the term $\B(\lambda)$ as
$$
\B(\lambda)=\B_\zs{1}(\lambda)+2\B_\zs{2}(\lambda)\,.
$$
The terms $\B_\zs{1}(\lambda)$ and $\B_\zs{2}(\lambda)$
are defined as
$$
\B_\zs{1}(\lambda)
=\frac{(b-a)^{2}}{d\sqrt{d}}\,
\sum^{d}_\zs{l=1}\,
\psi_\zs{1,l}(\lambda)\check{\eta}_\zs{l}
\quad\mbox{and}\quad
\B_\zs{2}(\lambda)
=
\frac{(b-a)^{2}}{d\sqrt{d}}\,
\sum^{d}_\zs{l=2}\,
\psi_\zs{2,l}(\lambda)\eta_\zs{l}\,,
$$
where
$$
\psi_\zs{1,l}(\lambda)=\sum^{d}_\zs{j=1}
\lambda(j)\phi^2_\zs{j}(z_\zs{l})
\quad\mbox{and}\quad
\psi_\zs{2,l}(\lambda)=\sum^{d}_\zs{j=1}
\lambda(j) \check{\upsilon}_\zs{j,l}
\,.
$$
So,
\begin{equation}
\label{sec:Prsm.11-nn}
\E_\zs{p,S}\,\B^2(\lambda)\le 
2 \E_\zs{p,S}\,\B^2_\zs{1}(\lambda)
+8
\E_\zs{p,S}\,\B^2_\zs{2}(\lambda)
\,.
\end{equation}
Taking into account property \eqref{sec:Ms.4nn_1}
and Bounyakovskii - Cauchy - Schwarz inequality we get
$$
\psi^2_\zs{1,l}(\lambda)\le 
\sum^{d}_\zs{j=1}\lambda^2(j)\phi^2_\zs{j}(z_\zs{l})
\sum^{d}_\zs{j=1}\phi^2_\zs{j}(z_\zs{l})
= \frac{d}{b-a}\,
\sum^{d}_\zs{j=1}\lambda^2(j)\phi^2_\zs{j}(z_\zs{l})\,.
$$
In view of properties
\eqref{sec:def_eta++_prs_0}
 we obtain that
\begin{align*}
\E_\zs{p,S}\,\B^{2}_\zs{1}(\lambda)
&=
\frac{(b-a)^{4}}{d^{3}}
\sum^{d}_\zs{l=1}\,\psi^2_\zs{1,l}(\lambda)\,\E_\zs{p,S}\,\check{\eta}^{2}_\zs{l}
\le 
\frac{(b-a)^{4}}{d^{3}}
\sum^{d}_\zs{l=1}\,\psi^2_\zs{1,l}(\lambda)\,\E_\zs{p,S}\,\eta^{4}_\zs{l}
\\[2mm]
&\le 
\frac{\sigma_\zs{1,*}\check{\m}(b-a)^{3}}{d^{2}}
\E_\zs{p,S}\,
\sum^{d}_\zs{j=1}\lambda^2(j)\sum^{d}_\zs{l=1}\,\sigma^{2}_\zs{l}\,\phi^2_\zs{j}(z_\zs{l})
\\[2mm]
&=\,\sigma_\zs{1,*}\,(b-a)\check{\m}\, \E_\zs{p,S}\,P_\zs{d}(\lambda)
\,.
\end{align*}
To estimate the last term in the right hand side of the inequality
\eqref{sec:Prsm.11-nn}, noting that the term $\psi_\zs{2,l}$ can be represented as
$$
\psi_\zs{2,l}(\lambda)=\sum^{l-1}_\zs{r=1}\,g_\zs{l,r}\,\eta_\zs{r}\,,
$$
where $g_\zs{l,r}=\sum^{d}_\zs{j=1}
\lambda(j)\,\phi_\zs{j}(z_\zs{l})\phi_\zs{j}(z_\zs{r})$, we use properties 
\eqref{sec:def_eta++_prs_0} to obtain
\begin{align*}
\E_\zs{p,S}\,\B^{2}_\zs{2}(\lambda)
&=
\frac{(b-a)^{4}}{d^{3}}
\sum^{d}_\zs{l=2}\,\E_\zs{p,S}\,\psi^2_\zs{2,l}(\lambda)\,\eta^{2}_\zs{l}
\le 
\frac{\sigma_\zs{1,*}(b-a)^{4}}{d^{3}}
\sum^{d}_\zs{l=2}\,\E_\zs{p,S}\psi^2_\zs{2,l}(\lambda)\,
\\[2mm]
&=
\frac{\sigma_\zs{1,*}(b-a)^{4}}{d^{3}}
\,
\sum^{d}_\zs{l=2}\,
\sum^{l-1}_\zs{r=1}\,
 g^{2}_\zs{l,r}\,\E_\zs{p,S}\sigma^{2}_\zs{r}\,
\\[2mm]
&\le
\frac{\sigma_\zs{1,*}(b-a)^{4}}{d^{3}}
\,
\sum^{d}_\zs{r=1}\,\E_\zs{p,S}\sigma^{2}_\zs{r}\,
\sum^{d}_\zs{l=1}\,
 g^{2}_\zs{l,r}
=\,\sigma_\zs{1,*}\,(b-a)\, \E_\zs{p,S}\,P_\zs{d}(\lambda)
\,.
\end{align*}
Hence Proposition \ref{Pr.sec:A.06-3}.
\fdem

Now we need the following moment bound.

\begin{proposition}\label{Pr.sec:App.1++0}
For any non random    $v_1,\ldots,v_\zs{d}$
\begin{equation}
\label{sec:Pr.00__+++1}
\E\,\left(\sum_{j=1}^{d} v_\zs{j}\,\eta_\zs{j,d}\right)^2
\,\le\,
\sigma_\zs{1,*}
\sum_{j=1}^{d}v_\zs{j}^2\,.
\end{equation}
\end{proposition}
\proof
Note that
\begin{align*}
\E\,\left(\sum_{j=1}^{d}\,v_\zs{j}\,\eta_\zs{j,d}\right)^2&=
\frac{b-a}{d}\,
\E\,\sum^{d}_\zs{l=1}\sigma^{2}_\zs{l}
\left(\sum^{d}_\zs{j=1}v_\zs{j} \phi_\zs{j}(z_\zs{l})\right)^{2}\\[2mm]
&\le 
\frac{\sigma_\zs{1,*}(b-a)}{d}\,
\sum^{d}_\zs{l=1}
\left(\sum^{d}_\zs{j=1}v_\zs{j} \phi_\zs{j}(z_\zs{l})\right)^{2}\,.
\end{align*}
By applying the orthonormal property \eqref{sec:Ms.4nn_1} we obtain the desired inequality. 
Hence Proposition~\ref{Pr.sec:App.1++0}.
\fdem

\medskip


\section{Proofs}\label{sec:proofs++}

\subsection{Proof of  Proposition \ref{Pr.sec:Prs.stp.times.1}}
First recall that
$$
\wh{S}_\zs{\iota_\zs{l}}=\frac{1}{A_\zs{\iota_\zs{l}}}\sum_{j=1}^{\iota_\zs{l}}Q_\zs{l,j}y_\zs{j-1}y_\zs{j}\,
\quad\mbox{and} \quad
\wt{S}_\zs{\iota_\zs{l}}=\min(\max(\wh{S}_\zs{\iota_\zs{l}},-1+\tilde{\epsilon}),1-\tilde{\epsilon})\,,
$$
where $\tilde{\epsilon}=1/(2+\ln n)$. Note that 
for sufficiently large $n$, for which we have $\tilde{\epsilon}<\epsilon$ and then $S(z_\zs{l})\in[-1+\tilde{\epsilon};1-\tilde{\epsilon}]$. We can write
\begin{align*}
|\wt{S}_\zs{\iota_\zs{l}}-S(z_\zs{l})|  \leq |\wh{S}_\zs{\iota_\zs{l}}-S(z_\zs{l})| 
 \leq  \frac{\displaystyle\sum_{j=k_\zs{1,l}}^{\iota_\zs{l}} y_\zs{j-1}^2 |S(x_\zs{j})-S(z_\zs{l})|}{\displaystyle\sum_{j=k_\zs{1,l}}^{\iota_\zs{l}}y_\zs{j-1}^2} 
  + 
  \vert 
  I_\zs{n} 
  \vert
  \,,
\end{align*}
where $I_\zs{n} = \displaystyle\sum_{j=k_\zs{1,l}}^{\iota_\zs{l}}\,y_\zs{j-1}\xi_\zs{j}/\displaystyle\sum_{j=k_\zs{1,l}}^{\iota_\zs{l}}y_\zs{j-1}^2$.
Taking into account that $\vert x_\zs{j}-z_\zs{l} \vert\le h$ for $k_{1,l}\le j\le k_{2,l}$, we obtain that for any 
$S\in\Theta_\zs{\epsilon,L}$,
$$
|\wh{S}_\zs{\iota_\zs{l}}-S(z_\zs{l})|\leq Lh + \vert I_\zs{n}\vert
\,.
$$
So,  for sufficiently large $n$
\begin{align}\nonumber
 \P_\zs{p,S}&\left(|\wt{S}_\zs{\iota_\zs{l}}-S(z_\zs{l})|>\epsilon_\zs{0}\right) 
 \leq  \P_\zs{p,S}\left(I_\zs{n}>\frac{\epsilon_\zs{0}}{2}\right) \\[2mm] \label{sec:UpperBnd_0011}
& \leq \P_\zs{p,S}\left(I_\zs{n}>\frac{\epsilon_\zs{0}}{2}\,,\,\Xi\right) +
 \P_\zs{p,S}\left(\Xi^{c} \right)\,,
\end{align}
where $\Xi=\left\{\left|\Upsilon_\zs{m_\zs{0},m_\zs{1}}(z_\zs{l})\right|\le  1/2\right\}$, $m_\zs{0}=k_\zs{1,l}-2$, $m_\zs{1}=\iota_\zs{l}-1$ and $\Upsilon_\zs{m_\zs{0},m_\zs{1}}(z_\zs{l})$ is defined in $\eqref{secPrsstptimesdef-0}$. Hence we obtain the following inequality on the set $\Xi$:
$$
\sum_{j=k_\zs{1,l}}^{\iota_\zs{l}} y_\zs{j-1}^2  = 
 (\iota_\zs{l}-k_\zs{1,l}+1) \left(\frac{1}{1-S^2(z_\zs{l})}+\Upsilon_\zs{m_\zs{0},m_\zs{1}}(z_\zs{l}) \right)
  \geq
\frac{ \q}{2}\,.
$$
Therefore, for any $\check{p}> 2$, 
\begin{align}\nonumber
 \P_\zs{p,S}\left(I_\zs{n}>\frac{\epsilon_\zs{0}}{2}\,,\,\Xi\right) 
& \leq  \P_\zs{p,S}\,\left(\left| \sum_{j=k_\zs{1,l}}^{\iota_\zs{l}} y_\zs{j-1}\xi_\zs{j}\right|>\frac{\q}{2}\right) \\[2mm]\label{upper-bd-010}
& \leq  \frac{2^{\check{p}}}{\q^{\check{p}}}\,
\E_\zs{p,S}\,\left|\sum_{j=k_\zs{1,l}}^{\iota_\zs{l}} y_\zs{j-1}\xi_\zs{j} \right|^{\check{p}}\,.
\end{align}

Using here the correlation inequality \eqref{sec:Bi.1++}
and the bound
\eqref{A.1},
we obtain that
$$
\max_\zs{1\le l\le d}
\sup_\zs{S\in \Theta_\zs{\e,L}}
\sup_\zs{p\in \cP}
\E_\zs{p,S}\,\left|\sum_{j=k_\zs{1,l}}^{\iota_\zs{l}} y_\zs{j-1}\xi_\zs{j} \right|^{\check{p}}
\le c_\zs{\check{p}}\,\q^{\check{p}/2}
\,.
$$
Applying this bound in
\eqref{sec:UpperBnd_0011}
and using Lemma \ref{Le.sec:Prs.stp.times.0-0}
we obtain  Proposition \ref{Pr.sec:Prs.stp.times.1}.
\endproof

\medskip

\medskip

\medskip 

\subsection{Proof of Proposition \ref{Pr.sec:Prs.stp.times.22}}
First note, that
$$
\P_\zs{p,S}(\Gamma^{c}) 
\le 
\sum^{d}_\zs{l=1}
\P_\zs{p,S}
\left(A_\zs{\iota_\zs{l},k_\zs{2,l}-1}<  H_\zs{l}\right)
\,.
$$
Moreover,  note that
in view of definition
 \eqref{secPrsstptimesdef-0} the term $A_\zs{\iota_\zs{l},k_\zs{2,l}-1}$ can be represented as
$$
A_\zs{\iota_\zs{l},k_\zs{2,l}-1}
=(m_\zs{1,l}-m_\zs{0,l})\,
\left(
\frac{1}{\gamma_\zs{l}}
+
\Upsilon_\zs{m_\zs{0,l},m_\zs{1,l}}(z_\zs{l})
\right)
\,,
$$
where $m_\zs{0,l}=\iota_\zs{l}-1$ and $m_\zs{1,l}=k_\zs{2,l}-2$.
Taking into account 
the definition of $H_\zs{l}$ in
\eqref{sec:Sp.8} and the fact
that $0<\wt{\gamma}_\zs{l},\gamma_\zs{l}\le 1$
and that $\vert\wt{\gamma}_\zs{l}-
\gamma_\zs{l}\vert \le 2\vert \wt{S}_\zs{\iota_\zs{l}}
-
S(z_\zs{l})
 \vert$, we obtain 
\begin{align*}
\P_\zs{p,S}
&\left(A_\zs{\iota_\zs{l},k_\zs{2,l}-1}<  H_\zs{l}\right)
= 
\P_\zs{p,S}
\left(
\frac{1}{\gamma_\zs{l}}
+
\Upsilon_\zs{m_\zs{0,l},m_\zs{1,l}}(z_\zs{l})
<\frac{1-\wt{\epsilon}}{\wt{\gamma}_\zs{l}}
\right)
\\[2mm]
&
\le
\P_\zs{p,S}
\left(
\left\vert
\frac{1}{\gamma_\zs{l}}
-
\frac{1}{\wt{\gamma}_\zs{l}}
\right\vert
>\frac{\wt{\epsilon}}{2}
\right)
+
\P_\zs{p,S}
\left(
\left\vert
\Upsilon_\zs{m_\zs{0,l},m_\zs{1,l}}(z_\zs{l})
\right\vert
>
\frac{\wt{\epsilon}}{2}
\right)\\[2mm]
&
\le
\P_\zs{p,S}
\left(
\left\vert
\wt{S}_\zs{\iota_\zs{l}}
-
S(z_\zs{l})
\right\vert
>\frac{\wt{\epsilon}^3}{4}
\right)
+
\P_\zs{p,S}
\left(
\left\vert
\Upsilon_\zs{m_\zs{0,l},m_\zs{1,l}}(z_\zs{l})
\right\vert
>
\frac{\wt{\epsilon}}{2}
\right)
\,.
\end{align*}
Applying here
Proposition \ref{Pr.sec:Prs.stp.times.1}
and Lemma \ref{Le.sec:Prs.stp.times.0-0}
we obtain
Proposition \ref{Pr.sec:Prs.stp.times.22}.
\endproof

\subsection{Proof of Proposition \ref{Pr.sec:Ms.0}}

Note that, for any $m\ge 1$
\begin{align*}
\E_\zs{p,S}\,\max_\zs{1\le j\le n}\, y^{4}_\zs{j}
&\le n^{\b/2}+
\sum^{n}_\zs{j=1}
\int^{+\infty}_\zs{n^{\b/2}}
\P_\zs{p,S}\left(y^{4}_\zs{j}\ge z\right)\d z\\[2mm]
&\le 
n^{\b/2}+n \max_\zs{1\le j\le n}\E_\zs{p,S}\,\vert y_\zs{j}\vert^{4m}
\int^{+\infty}_\zs{n^{\b/2}}
\, z^{-m}
\d z\\[2mm]
&=
n^{\b/2}+
\max_\zs{1\le j\le n}\E_\zs{p,S}\,\vert y_\zs{j}\vert^{4m}
\,
\frac{n^{1-\b(m-1)/2}}{m-1}
\,.
\end{align*}
Choosing here $m>1+2/\b$ and
using the bound  \eqref{A.1}
we obtain the property \eqref{sec:yyy_lim}. Hence Proposition \ref{Pr.sec:Ms.0}.
\endproof

\subsection{Proof of Theorem~\ref{Th.sec:Ms.1-0+++}}

First of all, note that on the set $\Gamma$ we can represent 
the empirical squared error $\Er_\zs{d}(\lambda)$ in the form
\begin{equation}\label{sec:Or.1}
\Er_\zs{d}(\lambda)=J_\zs{d}(\lambda)+
2\sum^{d}_\zs{j=1}\lambda(j)\theta'_\zs{j,d}+
\|S\|^2_\zs{d}-\delta\,
P_\zs{d}(\lambda)
\end{equation}
with $\theta'_\zs{j,d}=\wt{\theta}_\zs{j,d}-\theta_\zs{j,d}\wh{\theta}_\zs{j,d}$.
From \eqref{sec:Ms.5} we find that
$$
\theta'_\zs{j,d}=\theta_\zs{j,d}\zeta_\zs{j,d}+
\frac{b-a}{d}\wt{\eta}_\zs{j,d}
+
2\sqrt{\frac{b-a}{d}}\eta_\zs{j,d}\varpi_\zs{j,d}+
\varpi^2_\zs{j,d}\,,
$$
where $\wt{\eta}_\zs{j,d}=\eta^2_\zs{j,d}-\s_\zs{j,d}$. 
Now putting
\begin{equation}\label{sec:Or.2}
M(\lambda)=
\sum^{d}_\zs{j=1}\lambda(j)\,\theta_\zs{j,d}\,\zeta_\zs{j,d}\,,
\end{equation}
 we rewrite \eqref{sec:Or.1} as follows
\begin{align}\nonumber
\Er_\zs{d}(\lambda)
&=J_\zs{d}(\lambda)+
2M(\lambda)
+2
\frac{1}{\sqrt{d}}
\,\B(\lambda)
\\[2mm] \label{sec:Or.3}
&+2\Delta(\lambda)
+
\|S\|^2_\zs{d}-\delta\,
P_\zs{d}(\lambda)\,,
\end{align}
where $\B(\lambda)$ is given in \eqref{sec:Prsm.1},  $\Delta(\lambda)=\Delta_\zs{1}(\lambda)+\Delta_\zs{2}(\lambda)$,
$$
\Delta_\zs{1}(\lambda)=
\sum^{d}_\zs{j=1}\,\lambda(j)\,\varpi^2_\zs{j,d}
\quad\mbox{and}\quad
\Delta_\zs{2}(\lambda)=2
\sqrt{
\frac{b-a}{d}}
\sum^{d}_\zs{j=1}\,
\lambda(j)\,\eta_\zs{j,d}\varpi_\zs{j,d}\,.
$$
In view of Proposition \ref{Pr.sec:A.06-3},
 for any $\lambda\in\bbr^{d}$,
\begin{equation}\label{sec:Or.4}
\E_\zs{p,S}\,\Chi_\zs{\Gamma}\B^2(\lambda)
\le 
\,
10\sigma_\zs{1,*}\,(b-a)\check{\m} 
\,\E_\zs{p,S}\,
P_\zs{d}(\lambda)
\,.
\end{equation}
Note that the inequalities
\eqref{sec:seq_bound_sigma}
imply that
\begin{equation}\label{sec:Or.4++1}
P_\zs{0,d}(\lambda)
\le 
P_\zs{d}(\lambda)
\le 
P_\zs{1,d}(\lambda)
\,,
\end{equation}
where
$$
P_\zs{0,d}(\lambda)=\frac{\sigma_\zs{0,*}(b-a)\vert\lambda\vert^{2}}{d}
\quad\mbox{and}\quad
P_\zs{1,d}(\lambda)=\frac{\sigma_\zs{1,*}(b-a)\vert\lambda\vert^{2}}{d}
\,.
$$

For $\Delta_\zs{1}(\lambda)$, taking into account the properties of Fourier coefficients we obtain that
\begin{equation}\label{sec:Or.5}
\sup_\zs{\lambda\in [0,1]^{d}}|\Delta_\zs{1}(\lambda)|
\le
\sum^{d}_\zs{j=1}\varpi^2_\zs{j,d}=
 \|\varpi\|^2_\zs{d}\,.
\end{equation}
To estimate the term $\Delta_\zs{2}(\lambda)$ we recall that, for any $\varepsilon>0$ and  any $x,y\in\mathbb{R}$
\begin{equation}\label{sec:Or.5-1}
2xy\le \varepsilon x^2+\varepsilon^{-1}y^2\,.
\end{equation}
Therefore,
for some
$0<\varepsilon<1$, 
$$
|\Delta_\zs{2}(\lambda)|
\le 
\varepsilon 
\,\frac{b-a}{d}\,
 \sum^{d}_\zs{j=1}\,
\lambda^2(j)\, \eta^2_\zs{j,d}+
\frac{\|\varpi\|^2_\zs{d}}{\varepsilon}
= 
\varepsilon 
P_\zs{d}(\lambda)
+
\varepsilon \frac{\vert\B(\lambda^{2})\vert}{\sqrt{d}}
+
\frac{\|\varpi\|^2_\zs{d}}{\varepsilon}\,,
$$
where the vector $\lambda^2\in\Lambda_1$ as in \eqref{sec:Ms.9}.
Thus, for any $\lambda\in [0,1]^{d}$,
\begin{equation}\label{sec:Or.6}
|\Delta(\lambda)|\le 
\varepsilon 
P_\zs{d}(\lambda)
+
\varepsilon \frac{\B(\lambda^{2})\vert}{\sqrt{d}}+
2\varepsilon^{-1}\|\varpi\|^2_\zs{d}\,.
\end{equation}
Putting
$$
M_\zs{1}(\lambda)=2\frac{\B(\lambda)}{\sqrt{d}}+2\Delta(\lambda)\,,
$$
we can rewrite the empirical risk \eqref{sec:Or.3} as
\begin{equation}\label{sec:Or.6-1}
\Er_\zs{d}(\lambda)=J_\zs{d}(\lambda)+
2M(\lambda)
+M_\zs{1}(\lambda)
+
\|S\|^2_\zs{d}-\delta\,
P_\zs{d}(\lambda)\,.
\end{equation}
From \eqref{sec:Or.6} we obtain
$$
|M_\zs{1}(\lambda)|\le 
2\frac{\vert\B(\lambda)\vert}{\sqrt{d}}
+
2 
\frac{\vert\B(\lambda^{2})\vert}{\sqrt{d}}
+
2\varepsilon P_\zs{d}(\lambda)
+
4\varepsilon^{-1}
\|\varpi\|^2_\zs{d}\,.
$$
Moreover, setting
$$
\B^{*}=\sup_\zs{\lambda\in\Lambda}
\left(
\frac{\B^2(\lambda)}{ P_\zs{d}(\lambda)}
+
\frac{\B^2(\lambda^{2})}{ P_\zs{d}(\lambda^{2})}
\right)
$$
and taking into account that 
$P_\zs{d}(\lambda^2)\le P_\zs{d}(\lambda)$ for any $\lambda\in \Lambda$, 
we get 
$$
2\frac{\vert \B(\lambda)\vert}{\sqrt{d}}
+
2
\,\frac{\vert\B(\lambda^{2})\vert}{\sqrt{d}}
\le
2\varepsilon P_\zs{d}(\lambda)
+
\varepsilon^{-1}
\frac{\B^*}{d}
\,.
$$
By choosing $\varepsilon=\delta/4$ we find
\begin{equation}\label{sec:Or.7}
|M_\zs{1}(\lambda)|
\le 
\delta P_\zs{d}(\lambda)
+
\frac{16}{\delta}
\Upsilon_\zs{d}\,,
\quad
\Upsilon_\zs{d}
=
\frac{\B^*}{4 d}
+
\|\varpi\|^2_\zs{d}
\,.
\end{equation}
Now from \eqref{sec:Or.6-1}
we obtain that, for some
fixed $\lambda_0$ from $\Lambda$,
\begin{align*}
\Er_\zs{d}(\wh{\lambda})-\Er_\zs{d}(\lambda_0)\,&=\,J_\zs{d}(\wh{\lambda})\,-\,J_\zs{d}(\lambda_0)
+2\,M(\wh{\mu})
\\[2mm]
&
+M_\zs{1}(\wh{\lambda})
-\delta P_\zs{d}(\wh{\lambda})
-M_\zs{1}(\lambda_\zs{0})
+\delta P_\zs{d}(\lambda_\zs{0})\,,
\end{align*}
where $\wh{\mu}=\wh{\lambda}-\lambda_\zs{0}$.
By  the definition of $\wh{\lambda}$ in 
\eqref{sec:Ms.13} we obtain on the set $\Gamma$
\begin{align}\label{sec:Or.8}
\Er_\zs{d}(\wh{\lambda})&\le\Er_\zs{d}(\lambda_0)\,+\,
2\,M(\wh{\mu})+32\frac{\Upsilon_\zs{d}}{\delta}
+2\delta P_\zs{d}(\lambda_0)\,.
\end{align}
From
\eqref{sec:Or.4} and 
\eqref{sec:Or.4++1}
it follows that
\begin{align*}
\E_\zs{p,S}\,\Chi_\zs{\Gamma} \B^*&\le \sum_\zs{\lambda\in\Lambda}
\E_\zs{p,S}\,\Chi_\zs{\Gamma}\,
\left(
\frac{\B^2(\lambda)}{P_\zs{d}(\lambda)}
+
\frac{\B^2(\lambda^{2})}{P_\zs{d}(\lambda^{2})}
\right)\\[2mm]
&\le
10\sigma_\zs{1,*}(b-a)\check{\m}
\sum_\zs{\lambda\in\Lambda}
\left(
\frac{P_\zs{1,d}(\lambda)}{P_\zs{0,d}(\lambda)}
+
\frac{P_\zs{1,d}(\lambda^{2})}{P_\zs{0,d}(\lambda^{2})}
\right)\\[2mm]
&=
20\check{\m}(b-a)\nu\,\overline{\sigma}_\zs{*}
\,.
\end{align*}
and $\overline{\sigma}_\zs{*}
 =\displaystyle\frac{\sigma^{2}_\zs{1,*}}{\sigma_\zs{0,*}}\,.$ Therefore, for $0<\delta<1$ this inequality allows to bound $\Upsilon_\zs{d}$
as
\begin{equation}\label{sec:Or.9}
\E_\zs{p,S}\Chi_\zs{\Gamma}\Upsilon_\zs{d}
\le \frac{5\check{\m}(b-a)\overline{\sigma}_\zs{*}\nu}{d}
+
\u^{*}_\zs{d}\,,
\end{equation}
where $\u^{*}_\zs{d}$ is given by \eqref{sec:Ms.3-1}.

\noindent Now we study the second term on the right-hand side of inequality \eqref{sec:Or.8}.
For any  weight vector $\lambda \in\Lambda$
we set $\mu=\lambda-\lambda_\zs{0}$. Then we decompose this term as 
\begin{equation}\label{sec:Or.10}
M(\mu)=Z(\mu)+V(\mu)
\end{equation}
with
 $$
Z(\mu)=\sqrt{\frac{b-a}{d}}\,\sum^{d}_\zs{j=1}\mu(j)\,\theta_\zs{j,d}\eta_\zs{j,d}
\quad
\mbox{and}
\quad
V(\mu)=\sum^{d}_\zs{j=1}\,\mu(j)\,\theta_\zs{j,d}\varpi_\zs{j,d}\,.
$$
\noindent We define now the weighted discrete Fourier transformation of $S$, i.e.
we set
\begin{equation}\label{sec:Or.11}
\check{S}_\zs{\mu}=
\sum^{d}_\zs{j=1}\mu(j)\,\theta_\zs{j,d}\phi_\zs{j}\,.
\end{equation}
Now by using Propositioin~\ref{Pr.sec:App.1++0}
we can estimate the term $Z(\mu)$ as
\begin{equation}\label{sec:Or.12}
\E_\zs{p,S}\Chi_\zs{\Gamma}\,Z^2(\mu)\,
\le
\frac{\sigma_\zs{1,*}(b-a)}{d}
\,\|\check{S}_\zs{\mu}\|^2_\zs{d}
:=\sigma_\zs{1,*}(b-a)\D(\mu)\,.
\end{equation}
Moreover, by the inequalities \eqref{sec:Or.5-1} with 
$\varepsilon=\delta$ and \eqref{sec:Or.5} we can estimate $V(\mu)$ as follows
\begin{align}\label{sec:Or.13}
2V(\mu)=2
\sum^{d}_\zs{j=1}\,\mu(j)\,\theta_\zs{j,d}\varpi_\zs{j,d}
\le \delta\,\|\check{S}_\zs{\mu}\|^2_\zs{d}
+\frac{\|\varpi\|^2_\zs{d}}{\delta}\,.
\end{align}
Setting
$$
Z^*=\sup_\zs{\mu\in \Lambda-\lambda_\zs{0}}\frac{Z^2(\mu)}{\D(\mu)}\,,
$$
we obtain on the set $\Gamma$
\begin{equation}\label{sec:Or.14}
2M(\mu)\le\,2\delta \|\check{S}_\zs{\mu}\|^2_\zs{d}+
\frac{Z^*}{d\delta}+\frac{\|\varpi\|^2_\zs{d}}{\delta}\,.
\end{equation}
 Note now that from \eqref{sec:Or.12} it follows
 that
\begin{equation}\label{sec:Or.15}
\E_\zs{p,S}\,\Chi_\zs{\Gamma}\,Z^{*}
\le
\sum_\zs{\mu\in \Lambda-\lambda_\zs{0}}\frac{\E_\zs{p,S}\,\Chi_\zs{\Gamma}\,Z^2(\mu)}{\D(\mu)}
\le \nu\sigma_\zs{1,*}(b-a)\,.
\end{equation}
Now we estimate the first term on the right-hand side of the inequality \eqref{sec:Or.14}. 
On the set $\Gamma$ we have
\begin{align}\nonumber
\|\check{S}_\zs{\mu}\|^2_\zs{d}-
\|\wh{S}_\zs{\mu}\|^2_\zs{d}&=
\sum^{d}_\zs{j=1}\,\mu^2(j)
(\theta^2_\zs{j,d}-\wh{\theta}^2_\zs{j,d})\,
\le 
-\,2\sum^{d}_\zs{j=1}\,\mu^2(j)\,
\theta_\zs{j,d}\,\zeta_\zs{j,d}\\[2mm] \label{sec:Or.16}
&=-2 Z_\zs{1}(\mu)-2 V_\zs{1}(\mu)\,,
\end{align}
where
$$
Z_\zs{1}(\mu)=\sqrt{\frac{b-a}{d}}
\sum^{d}_\zs{j=1}\,\mu^2(j)
\theta_\zs{j,d}\eta_\zs{j,d}
\quad
\mbox{and}
\quad
V_\zs{1}(\mu)=\sum^{d}_\zs{j=1}\,\mu^2(j)\,\theta_\zs{j,d}
\varpi_\zs{j,d}\,.
$$
Taking into account that $|\mu(j)|\le 1$,
similarly to inequality \eqref{sec:Or.12}, we find 
$$
\E_\zs{p,S}\,\Chi_\zs{\Gamma}\,Z^2_\zs{1}(\mu)\,\le\,\sigma_\zs{1,*}\,\D(\mu)\,.
$$
Moreover, for the random variable
$$
Z^*_\zs{1}=\sup_\zs{\mu\in \Lambda-\lambda_\zs{0}}\frac{Z^2_\zs{1}(\mu)}{\D(\mu)}\,,
$$
we obtain the same upper bound as in \eqref{sec:Or.15}, i.e. 
\begin{equation}\label{sec:Or.17}
\E_\zs{p,S}\,Z^*_\zs{1}\,\Chi_\zs{\Gamma}\,\le \nu\sigma_\zs{1,*}(b-a)\,.
\end{equation}
Furthermore, similarly to \eqref{sec:Or.13} we estimate the second term in \eqref{sec:Or.16} as
$$
2\vert V_\zs{1}(\mu)\vert\le \delta\|\check{S}_\zs{\mu}\|^2_\zs{d}+
\frac{\|\varpi\|^2_\zs{d}}{\delta}\,.
$$
Therefore, on the set $\Gamma$ 
\begin{align*}
\|\check{S}_\zs{\mu}\|^2_\zs{d}\le
\|\wh{S}_\zs{\mu}\|^2_\zs{d}\,
+2\delta \|\check{S}_\zs{\mu}\|^2_\zs{d}+
\frac{Z^*_\zs{1}}{d\delta}
+
\frac{\|\varpi\|^2_\zs{d}}{\delta}\,,
\end{align*}
i.e.
\begin{align}\label{sec:Or.18}
\|\check{S}_\zs{\mu}\|^2_\zs{d}\le
\frac{1}{1-2\delta}\,
\|\wh{S}_\zs{\mu}\|^2_\zs{d}\,
+
\frac{1}{(1-2\delta)\delta}
\left(
\frac{Z^*_\zs{1}}{d}
+
\|\varpi\|^2_\zs{d}
\right)\,.
\end{align}
Using this inequality in \eqref{sec:Or.14} and putting 
$Z^*_\zs{2}=Z^*+Z^*_\zs{1}$ yield on the set $\Gamma$
\begin{align*}
2M(\wh{\mu})&\le\,\frac{2\delta}{1-2\delta }
\|\wh{S}_\zs{\wh{\mu}}\|^2_\zs{d}
+
\frac{1}{\delta(1-2\delta)}
\left(\frac{Z^*_\zs{2}}{d}+
\|\varpi\|^2_\zs{d}
\right)\\[2mm]
&\le\,\frac{4\delta(\Er_\zs{d}(\wh{\lambda})+\Er_\zs{d}(\lambda_0))}{1-2\delta}
+
\frac{1}{\delta(1-2\delta)}
\left(
\frac{Z^*_\zs{2}}{d}+
\|\varpi\|^2_\zs{d}
\right)\,.
\end{align*}
Therefore from the preceding inequality and \eqref{sec:Or.8} we obtain
\begin{align*}
\Er_\zs{d}(\wh{\lambda})\Chi_\zs{\Gamma}&
\le \, \frac{1+2\delta}{1-6\delta }
\Er_\zs{d}(\lambda_0)\Chi_\zs{\Gamma}\,\,+\,
\frac{32(1-2\delta) }{\delta(1-6\delta)}\Upsilon_\zs{d}\Chi_\zs{\Gamma}\\[2mm]
&
+
\frac{1}{\delta(1-6\delta)}
\left(\frac{Z^*_\zs{2}}{d}+
\|\varpi\|^2_\zs{d}
\right)
\Chi_\zs{\Gamma}
+\frac{2\delta(1-2\delta)}{1-6\delta}\,P_\zs{d}(\lambda_\zs{0})\Chi_\zs{\Gamma}
\end{align*}
and through the  inequalities \eqref{sec:Or.9},
 \eqref{sec:Or.15} and \eqref{sec:Or.17} we estimate the empirical risk
as
\begin{align*}
\E_\zs{p,S}\Er_\zs{d}(\wh{\lambda})\Chi_\zs{\Gamma}&
\le \, \frac{1+2\delta}{1-6\delta}
\E_\zs{p,S}\Er_\zs{d}(\lambda_\zs{0})\Chi_\zs{\Gamma}\,+\,
\frac{32(1-2\delta)}{\delta(1-6\delta)}\,\left(
\frac{5\check{\m}\overline{\sigma}_\zs{*}\nu(b-a)}{d}+\u^{*}_\zs{d} 
\right)\\[2mm]
&
+\frac{1}{\delta(1-6\delta)}
\left(
\frac{2\nu\sigma_\zs{1,*}(b-a)}{d}+\u^{*}_\zs{d}
\right)
+\frac{2\delta(1-2\delta)}{1-6\delta}\,
\E_\zs{p,S}\Chi_\zs{\Gamma} P_\zs{d}(\lambda_\zs{0})\,. 
\end{align*}
Taking into account that $\sigma_\zs{*,1}\le \overline{\sigma}_\zs{*}$ and that $1-6\delta>1/2$ for $0<\delta<1/12$, we get 
\begin{align*}
\E_\zs{p,S}\Er_\zs{d}(\wh{\lambda})\Chi_\zs{\Gamma}&
\le \, \frac{1+2\delta}{1-6\delta}
\E_\zs{p,S}\Er_\zs{d}(\lambda_\zs{0})\Chi_\zs{\Gamma}+
\frac{320}{\delta}
\left(
\frac{(\check{\m}+1)\overline{\sigma}_\zs{*}\nu(b-a)}{d}
+ \u^{*}_\zs{d}
\right)\\[2mm]
&+\frac{2\delta(1-2\delta)}{1-6\delta}\,
\E_\zs{p,S}\Chi_\zs{\Gamma} P_\zs{d}(\lambda_\zs{0})\,. 
\end{align*}
By applying Lemma~\ref{Le.sec:App.2} with $\varepsilon=2\delta$
we get that
\begin{align*}
\E_\zs{p,S}\Er_\zs{d}(\wh{\lambda})\Chi_\zs{\Gamma}
&
\le \, \frac{1+4\delta}{1-6\delta}
\E_\zs{p,S}\Er_\zs{d}(\lambda_\zs{0})\Chi_\zs{\Gamma}
+
\frac{320}{\delta}
\left(
\frac{(\check{\m}+1)\overline{\sigma}_\zs{*}\nu(b-a)}{d}
+3\u^{*}_\zs{d}
\right)
\\[2mm]&
+
10\delta\sqrt{\sigma_\zs{1,*}\check{\m}\P_\zs{p,S}(\Gamma^c)}\,.
\end{align*}
Taking into account the definition of $\check{\m}$ in
\eqref{sec:def_eta++_prs_0} and that $\m^{*}_\zs{4}\le 3\varsigma^{2}$, then by
replacing 
$$
\E_\zs{p,S}\Er_\zs{d}(\wh{\lambda})\Chi_\zs{\Gamma}
\quad\mbox{and}\quad
\E_\zs{p,S}\Er_\zs{d}(\lambda_\zs{0})\Chi_\zs{\Gamma}
$$
by 
$$
\E_\zs{p,S}\|\wh{S}_\zs{*}-S\|^2_\zs{d}-\|S\|^2_\zs{d}\P_\zs{p,S}(\Gamma^c)
\quad\mbox{and}\quad
\E_\zs{p,S}\|\wh{S}_\zs{\lambda_\zs{0}}-S\|^2_\zs{d}-\|S\|^2_\zs{d}\P_\zs{p,S}(\Gamma^c)
$$
respectively,
we come to the inequality \eqref{sec:Mrs.1--00-22}. Hence 
Theorem~\ref{Th.sec:Ms.1-0+++}.
\endproof

\bigskip

\medskip

{\bf Acknowledgements.}
The last author was  partially supported by
by 
  the Russian Federal Professor program (project no. 1.472.2016/1.4, the  Ministry of Education and Science of the Russian Federation),
 the research project no. 2.3208.2017/4.6
(the  Ministry of Education and Science of the Russian Federation), by  RFBR Grant 16-01-00121 A and
 by "The Tomsk State University competitiveness improvement program" grant 8.1.18.2018.

\medskip

\medskip

\setcounter{section}{0}
\renewcommand{\thesection}{\Alph{section}}

\section{Appendix}\label{sec:A}

\subsection{Burkh\"older inequality}

\noindent
We need the following  from \cite{Shiryaev2004}.

\begin{proposition}\label{Pr.A.2-0}
Let $(M_\zs{k})_\zs{1\le k\le n}$ be a martingale. Then 
for any $q>1$
\begin{equation}\label{A.1-1}
\E\,|M_\zs{n}|^{q}\le \b^{*}_\zs{q}
 \E \left(\sum^{n}_\zs{j=1}
(M_\zs{j}-M_\zs{j-1})^{2}\right)^{q/2}\,,
\end{equation}
where the coefficient $\b^{*}_\zs{q} = 18(q)^{3/2}/(q-1)^{1/2}.$ 
\end{proposition}

\medskip

\subsection{Properties of the sequential procedures}

\begin{lemma}\label{Le.sec:App.eta_prps}
The properties \eqref{sec:def_eta++_prs_0} hold for the random variables $(\eta_\zs{l})_\zs{1\le l\le d}$
defined in
\eqref{sec:def_eta}.
\end{lemma}
\proof
First, we set $\cF_\zs{j}=\sigma\{\xi_\zs{1},\ldots,\xi_\zs{j}\}$ for $1\le j\le n$ and as usual \\
$\cF_\zs{0}=\{\Omega,\emptyset\}$.
Moreover, note that 
$$
\eta_\zs{l}=\sum^{n}_\zs{j=1}\,\check{t}_\zs{l,j}\,\xi_\zs{j}
\quad\mbox{and}\quad
\check{t}_\zs{l,j}=\sigma^{2}_\zs{l}\left(\Chi_\zs{\{\iota_\zs{l}\le j<\check{\tau}_\zs{l}\}}\check{Q}_\zs{l,j}+
 \Chi_\zs{\{j=\check{\tau}_\zs{l}\}}\,\check{\varkappa}_\zs{l}\,
 \check{Q}_\zs{l,\check{\tau}_\zs{l}}\right)
 \,.
$$ 
Taking into account that $\check{t}_\zs{l,j}$ is $\cF_\zs{j-1}$ - measurable for any $1\le j\le n$ and
$$
\sum^{n}_\zs{j=1}\check{t}^{2}_\zs{l,j}=\sigma^{2}_\zs{l}\,.
$$
Note also that $\cG_\zs{l}=\sigma\{\eta_\zs{1},\ldots,\eta_\zs{l-1},\sigma_\zs{l},\}\subset\cF_\zs{\iota_\zs{l}}$. Noting that
$$
\E\,\left(\eta_\zs{l}\vert \cF_\zs{\iota_\zs{l}}\right)=0
\quad\mbox{and}\quad
\E\,\left(\eta^{2}_\zs{l}\vert \cF_\zs{\iota_\zs{l}}\right)=1\,,
$$
we obtain the first two equalities in
 \eqref{sec:def_eta++_prs_0}.
As to the last inequality, note that
 through \eqref{A.1-1} we can write
 $$
 \E_\zs{p,S}\left(\left(\sum^{n}_\zs{j=1}\,\check{t}_\zs{l,j}\,\xi_\zs{j}\right)^{4} \vert \cF_\zs{\iota_\zs{l}}\right)
 \le 
 \b^{*}_\zs{4}\,
 \E_\zs{p,S}\left(\left(\sum^{n}_\zs{j=1}\,\check{t}^{2}_\zs{l,j}\,\xi^{2}_\zs{j}\right)^{2} \vert \cF_\zs{\iota_\zs{l}}\right)
 \,.
 $$
Now, note that
$$
\left(\sum^{n}_\zs{j=1}\,\check{t}^{2}_\zs{l,j}\,\xi^{2}_\zs{j}\right)^{2}\le 2
\sigma^{4}_\zs{l}
+
2
\left(\sum^{n}_\zs{j=1}\,\check{t}^{2}_\zs{l,j}\,\wt{\xi}_\zs{j}\right)^{2}
$$
where $\wt{\xi}_\zs{j}=\xi^{2}_\zs{j}-1$. Taking into account that
$$
\E_\zs{p,S}\left(
\left(\sum^{n}_\zs{j=1}\,\check{t}^{2}_\zs{l,j}\,\wt{\xi}_\zs{j}\right)^{2}
\vert \cF_\zs{\iota_\zs{l}}\right)
=\E_\zs{p}\,\wt{\xi}^{2}_\zs{1}
\,
\sum^{n}_\zs{j=1}\,\check{t}^{4}_\zs{l,j}\,
\le 
\sigma^{4}_\zs{l}\,\E_\zs{p}\,\wt{\xi}^{2}_\zs{1}\,,
$$
we obtain the last bound in \eqref{sec:def_eta++_prs_0}. Hence Lemma \ref{Le.sec:App.eta_prps}.

\fdem

\subsection{Correlation inequality}

Now we give the correlation inequality from \cite{GaltchoukPergamenshchikov2013}.

\begin{proposition}\label{Pr.sec:corr_ineq}
 Let $(\Omega,\cF,(\cF_\zs{j})_\zs{1\le j\le n},\P)$
be a filtered probability space and
$(X_\zs{j}, \cF_\zs{j})_\zs{1\le j\le n}$ be a sequence of random
variables such that for some $p\ge 2$
$$
\max_\zs{1\le j\le n}\,\E\,|X_\zs{j}|^{p}\,<\,\infty\,.
$$
Define
$$
b_\zs{j,n}(p)=
\left(
\E\,
(
|X_\zs{j}|\,\sum^{n}_\zs{k=j}
|\E\,(X_\zs{k}|\cF_\zs{j})|
)^{p/2}
\right)^{2/p}\,.
$$
 Then
\begin{equation}\label{sec:Bi.1++}
\E\,|
\sum^{n}_\zs{j=1}\,X_\zs{j}
\,|^{p}
\le\,(2p)^{p/2}
\left(
\sum^{n}_\zs{j=1}\,
b_\zs{j,n}(p)
\right)^{p/2}\,.
\end{equation}
\end{proposition}

\medskip
\subsection{Upper bound for the penalty term}\label{subsec:App.3}
\begin{lemma}\label{Le.sec:App.2}
For sufficiently large $n$ and $0<\varepsilon<1$,
\begin{align*}
\E_\zs{p,S} P_\zs{d}(\lambda)
&\le 
\,\frac{1}{1-\varepsilon}\E_\zs{p,S}\,\Er_\zs{d}(\lambda)\Chi_\zs{\Gamma}+
\frac{\u^*_\zs{d}}{(1-\varepsilon)\varepsilon }\\[2mm]
&
+\frac{10}{1-\varepsilon}\,\sqrt{\sigma_\zs{1,*}\check{\m}\P_\zs{p,S}(\Gamma^c)}\,.
\end{align*}
\end{lemma}
\proof
Indeed, by the definition of $\Er_\zs{d}(\lambda)$ on the set $\Gamma$
 we have
\begin{align*}
\Er_\zs{d}(\lambda)\,&=\,\sum^{d}_\zs{j=1}
\left((\lambda(j)-1)\theta_\zs{j,d}+\lambda(j)\zeta_\zs{j,d}\right)^2\\[2mm]
&=\sum^{d}_\zs{j=1}\,
\left((\lambda(j)-1)\theta_\zs{j,d}+\lambda(j)\varpi_\zs{j,d}+\lambda(j)
\sqrt{\frac{b-a}{d}}
\eta_\zs{j,d}\right)^2\,.
\end{align*}
Therefore, putting
$$
I_\zs{1}=\sum^{d}_\zs{j=1}
\lambda(j)(\lambda(j)-1)\theta_\zs{j,d}\eta_\zs{j,d}
\quad
\mbox{and}
\quad
I_\zs{2}=\sum^{d}_\zs{j=1}\lambda^2(j)\varpi_\zs{j,d}\eta_\zs{j,d}
\,,
$$
we get on the set $\Gamma$ the following lower bound for the empirical risk
$$
\Er_\zs{d}(\lambda)
\ge
\,\frac{b-a}{d}\,
\sum^{d}_\zs{j=1}\lambda^2(j)\,\eta^2_\zs{j,d}
+2\sqrt{\frac{b-a}{d}}\,I_\zs{1}+2\sqrt{\frac{b-a}{d}}\,I_\zs{2}\,.
$$

Taking into account that for $0<\varepsilon<1$,
$$
2\sqrt{\frac{b-a}{d}}
|I_\zs{2}|\le 
\,\varepsilon\,\frac{b-a}{d}\sum^{d}_\zs{j=1}\lambda^2(j)\eta^2_\zs{j,d}
+\frac{\|\varpi\|^2_\zs{d}}{\varepsilon}\,,
$$
 we get 
\begin{equation}\label{sec:App.1}
\Er_\zs{d}(\lambda)\,\ge\,(1-\varepsilon)
\,\frac{b-a}{d}
\sum^{d}_\zs{j=1}\lambda^2(j)\,\eta^2_\zs{j,d}
+2\sqrt{\frac{b-a}{d}}\,I_\zs{1}-
\frac{\|\varpi\|^2_\zs{d}}{\varepsilon}\,.
\end{equation}

Let us consider the first term in \eqref{sec:App.1}, then we have
\begin{align*}
\E_\zs{p,S}\Chi_\zs{\Gamma}\sum^{d}_\zs{j=1}\lambda^2(j)\,\eta^2_\zs{j,d}=
\E_\zs{p,S}\sum^{d}_\zs{j=1}\lambda^2(j)\,\s_\zs{j,d}-
\E_\zs{p,S}\Chi_\zs{\Gamma^c}\sum^{d}_\zs{j=1}\lambda^2(j)\,\eta^2_\zs{j,d}\,.
\end{align*}

Using the correlation inequality \eqref{sec:Bi.1++}
and the upper bound for the fourth moment in 
\eqref{sec:def_eta++_prs_0}
we obtain
\begin{equation}\label{sec:App.2--00}
\E_\zs{p,S}\eta^4_\zs{j,d}\le 64 \check{\m}\,\sigma^2_\zs{1,*}
\,.
\end{equation}
This implies
\begin{equation}\label{sec:App.3}
\E_\zs{p,S}\Chi_\zs{\Gamma}\sum^{d}_\zs{j=1}\lambda^2(j)\,\eta^2_\zs{j,d}
\ge
\E_\zs{p,S}\sum^{d}_\zs{j=1}\lambda^2(j)\,\s_\zs{j,d}
-
8\sigma_\zs{1,*}\,d\,\sqrt{\check{\m}\,\P_\zs{p,S}(\Gamma^c)}\,.
\end{equation}

Therefore, we obtain
$$
\frac{b-a}{d}
\E_\zs{p,S}\Chi_\zs{\Gamma}\sum^{d}_\zs{j=1}\lambda^2(j)\,\eta^2_\zs{j,d}
\ge
\E_\zs{p,S}\, P_\zs{d}(\lambda)
-
8(b-a)
\sigma_\zs{1,*}\,\sqrt{\check{\m}\P_\zs{p,S}(\Gamma^c)}
\,.
$$
Moreover, taking into account that $\E_\zs{p,S}I_\zs{1}=0$ and 
in view of Proposition \eqref{Pr.sec:App.1++0}
$$
\E_\zs{p,S}I^{2}_\zs{1}\le 
\sigma_\zs{1,*}\,\|S\|^{2}_\zs{d}
\,.
$$
So, recalling that 
that $\|S\|_\zs{d} \le b-a$, we estimate 
$\E_\zs{p,S}I_\zs{1}\Chi_\zs{\Gamma}$ as
\begin{align*}
\left|\E_\zs{p,S}\,I_\zs{1}\Chi_\zs{\Gamma}\right|
=
\left|\E_\zs{p,S}\,I_\zs{1}\Chi_\zs{\Gamma^c}\right|
\le \,\sqrt{\sigma_\zs{1,*}}
\sqrt{\P_\zs{p,S}(\Gamma^c)}\,.
\end{align*}
\noindent Hence Lemma~\ref{Le.sec:App.2}.
\endproof

\medskip
\subsection{Properties of the model \eqref{sec:In.1}}

\begin{lemma}\label{Le.A.1}
For all $t\in \mathbb{N}^*$ and ~  $0<\e<1,$ the random variables $y_k$ in \eqref{sec:In.1}
 satisfy the following :
\begin{equation}\label{A.1}
\sup_\zs{n\ge 1}\,\sup_\zs{0\le k\le n}\,
\sup_\zs{S\in\Theta_\zs{\e,L}}\,
\E_\zs{p,S}\,y^{2t}_\zs{k}\, < \infty.
\end{equation}
\end{lemma}
\noindent
\proof This lemma is shown in \cite{Arkoun2011} (Lemma A.1).
\endproof

\noindent
We set
\begin{equation}
\label{secPrsstptimesdef-0}
\Upsilon_\zs{m_\zs{0},m_\zs{1}}(z_\zs{l})
=
\frac{1}{m_\zs{1}-m_\zs{0}}
\sum^{m_\zs{1}}_\zs{j=m_\zs{0}+1}
y^{2}_\zs{j}
-
\frac{1}{\gamma_\zs{l}}
\,,
\end{equation}
where $(k_\zs{1,l}-2)_\zs{+}\le m_\zs{0}<m_\zs{1}\le k_\zs{2,l}$, $(a)_\zs{+}$ is positive part of $a$ and $\gamma_\zs{l}$ is defined in \eqref{sec:Sp.6}.

\begin{lemma}\label{Le.sec:Prs.stp.times.0-0}
Assume that the bounds $m_\zs{0}$ and $m_\zs{1}$ in \eqref{secPrsstptimesdef-0} 
are such that  for some $0<\epsilon_\zs{1}<1/2$
$$
\liminf_\zs{n\to\infty}n^{\epsilon_\zs{1}}(m_\zs{1}-m_\zs{0})>0
\,.
$$
Then, for any $\b>0$ 
\begin{equation}\label{sec:Prs.stp.times.0-0}
\lim_\zs{n\to\infty}\,n^{\b}
\max_\zs{1\le l\le d}
\,
\sup_\zs{S\in \Theta_\zs{\e,L}}
\sup_\zs{p\in \cP}
\P_\zs{p,S}\left(
\left\vert
\Upsilon_\zs{m_\zs{0},m_\zs{1}}(z_\zs{l})
\right\vert
>\epsilon_\zs{0}\right)
=0\,,
\end{equation}
where $\epsilon_\zs{0}=\epsilon_\zs{0}(n)\to 0$ as $n\to\infty$ is such that 
$\lim_\zs{n\to\infty}n^{\check{\delta}}\epsilon_\zs{0}=\infty$
for any $\check{\delta}>0$.
\end{lemma}
\noindent
\proof
This lemma is shown in   \cite{Arkoun2011} (Lemma A.2).
\endproof

\medskip
\subsection{Properties of the norms}

\begin{lemma}\label{Le.sec:A.1-1}
Let $f$ be an absolutely continuous $[a,b]\to\bbr$ function with
$\|\dot{f}\|<\infty$ and $g$ be a simple $[a,b]\to\bbr$ function of the form
$$
g(t)=\sum_{j=1}^p\,c_\zs{j}\,\chi_{(t_{j-1},t_j]}(t),$$
where $c_\zs{j}$ are some constants. Then 
for all $ \wt{\varepsilon}>0$,  the function $\Delta=f-g$
satisfies the following inequalities
$$
\|\Delta\|^{2}\le (1+\wt{\varepsilon})\|\Delta\|^{2}_\zs{d}
+
\left(1+\frac{1}{\wt{\varepsilon}}\right)\frac{\|\dot{f}\|^{2}}{d^{2}}(b-a)^2\,,
\quad
$$
and
$$
\|\Delta\|^{2}_\zs{d}\le (1+\wt{\varepsilon})\|\Delta\|^{2}
+
\left(1+\frac{1}{\wt{\varepsilon}}\right)\frac{\|\dot{f}\|^{2}}{d^{2}}(b-a)^2\,.
$$
\end{lemma}
\noindent 
\proof
Lemma
\ref{Le.sec:A.1-1}
is proven in 
\cite{KonevPergamenshchikov2015}.
(Lemma A.2.)
\endproof

\medskip

\medskip

\bibliographystyle{plain}


\end{document}